\magnification=1200
\hsize=6truein
\hoffset=.5truein
\vsize=8 true in
\baselineskip=12pt plus 2pt minus 2pt
\parskip=3pt
\tolerance= 10000
\font\rmtwelve=cmbx10 at 12pt

{\rmtwelve Uniqueness of symplectic canonical class,
surface cone and symplectic cone of $4-$manifolds with $b^+=1$}
\bigskip
\bigskip
\centerline {Tian-Jun Li \& Ai-Ko Liu}
\bigskip
\bigskip

\noindent{\bf Abstract}. Let $M$ be a closed oriented smooth
$4-$manifold admitting symplectic structures.
If $M$ is minimal and has $b^+=1$, we prove that
there is a unique symplectic
canonical class up to sign, and any real second cohomology
class of positive square is represented by symplectic forms.
Similar results hold when $M$ is not minimal.
\bigskip
\bigskip

\noindent{\bf \S1. Introduction}

\medskip

Let $M$ be a smooth, closed oriented $4-$manifold. An orientation-compatible
symplectic structure
on $M$ is a closed 2-form $\omega$ such that $\omega\wedge \omega$
is nowhere vanishing and agrees with  the orientation.
Let $\Omega_M$ be the moduli space of such 2-forms.
In the first part of this paper, we devote ourselves
to the understanding of the topology of this moduli space,
which   can be studied in three ways.

First of all, on $\Omega_M$, there is a natural equivalence relation, the
deformation equivalence.
 $\omega_1$ and $\omega_2$ in $\Omega_M$
are said to be deformation equivalent  if there is an
orientation-preserving  diffeomorphism
$\phi$ such that $\phi^*\omega_1$ and $\omega_2$ are connected by
a path of symplectic forms.  Clearly, the group of orientation-preserving
diffeomorphisms
on the set of connected components of $\Omega_M$, and
the number of
deformation classes of symplectic structures is just the number of
the orbits of  this action.

Secondly, there is a map of canonical class $K: \Omega_M\longrightarrow
H^2(M;{\bf Z})$.
Any symplectic structure
determines a homotopy class of compatible almost complex structures
on the cotangent bundle,
whose first Chern class is called the symplectic canonical class.
For each symplectic canonical class $K$, if
 we let $\Omega_{M,K}$ be the subset of
$\Omega_M$,
whose elements have  $K$ as the symplectic canonical class,
then ${\Omega_M}$ is the disjoint union of
the $\Omega_{M,K}$. There is also a natural equivalence
relation on the set of symplectic canonical classes.
 We say two symplectic canonical classes $K_1$ and $K_2$ are equivalent
if there is an orientation-preserving
 diffeomorphism $\phi$ such that $\phi^* K_1=\pm K_2$.
Symplectic structures in a connected component have the same
symplectic canonical class.
Moreover  if two symplectic structures are
related by an orientation-preserving
 diffeomorphism, so are their symplectic canonical
classes. Therefore the set of deformation equivalence classes of
orientation-compatible symplectic structures maps onto the set of
equivalent classes of symplectic canonical classes, and can be
understood via the latter.

Thirdly, by taking the cohomology class, we have a
projection $CC:\Omega_M\longrightarrow H^2(M;{\bf R})$. The image of
this projection is called the symplectic cone of $M$, and is denoted
by ${\cal C}_M$.
For each symplectic canonical class,
if we define the $K-$symplectic cone
$${\cal C}_{M,K}=\{e\in {\cal C}_M|
 \hbox{ $e=[\omega]$ for some $\omega \in \Omega_K $ }\}.$$
then ${\cal C}_M$ is just the union of
the ${\cal C}_{M,K}$
(this union is in fact a disjoint union, see Proposition 4.1).

It turns out that the determination of the
number of deformation
classes of symplectic structures is a very hard  problem.
 The first breakthrough was made by Taubes [T3],  who showed that there is
one deformation class on $CP^2$. We [LL2] later
showed that the uniqueness holds for all $S^2-$bundles and their blow ups.
However, for other 4-manifolds, it is not even known whether the
number of deformation classes is always finite.

Here we are able to provide some evidence for the finiteness by
showing that the the image of the map $K$ is always finite.
When $b^+(M)>1$ ($b^+(M)$ is the dimension of a maximal positive subspace of
$H^2(M;{\bf R})$), it follows from [T1] and [W] that
there are finitely many symplectic canonical classes.
In addition, any 4-manifolds with K\"ahler structures  has only one equivalent class
of symplectic canonical classes (see [B] and [FM2]).
More recently, examples of 4-manifolds with $b^+>1$ and Inequivalent classes of symplectic canonical
classes were first obtained in [MT], and later in [Le],  [Sm] and [V] (see also [C], [KK] and [Ma] for the
recent results for the moduli space of complex structures on 4-manifolds).
In this paper, our first main result  completely
settles this issue for the case $b^+=1$.

\noindent{\bf Theorem 1}. Let $M$ be a smooth, closed oriented
$4-$manifold with $b^+=1$ and suppose $\Omega_M$ is not empty.
Then the number of equivalent classes of  the symplectic canonical classes
is one.
Furthermore, if $M$ is minimal, there is a unique symplectic
canonical class up to sign.

Here $M$ is said to be
(smoothly) minimal if ${\cal E}_M$
is the empty set,
where ${\cal E}_M$ is the set of cohomology classes whose Poincar\'e duals are
represented by
smoothly embedded spheres of squares $-1$. Consequently we obtain

\noindent{\bf Corollary 1}. Let $M$ be a smooth, closed oriented
$4-$manifold. The number of equivalent classes of the symplectic canonical
classes is finite.

To the contrary, in higher dimensions,
there can be infinite number of equivalent classes.
Let us briefly indicate the construction of such a manifold. Using their knot
surgery,
Fintushel and Stern [FS2] have constructed infinite number of   symplectic $4-$manifolds homeomorphic to $K3$, whose diffeomorphism types are
distinguished by the number of their Seiberg-Witten basic classes.
As can be shown with the methods in [R1], the products
of those manifolds with $S^2$ are diffeomorphic, and it  carries
infinite number of equivalent classes.

We also have some substantial results about the
map $CC$. More precisely, for a
4-manifold with $b^+=1$ and non-empty symplectic cone, we can
characterize the $K-$symplectic cones ${\cal C}_{M,K}$,
and from which, obtaining the
characterization of the symplectic cone ${\cal C}_M$ itself.

To determine the $K-$symplectic cone, the first
step is to investigate the set of `$K-$stable'
classes of symplectic surfaces, $A_{M,K}$, which we define now.
Let $A_{M,\omega}$ be the set of $s\in H^2(M;{\bf Z})$
whose Poincar\'e duals can
be represented by embedded $\omega-$symplectic surfaces.  Define
$$A_{M,K}=\{e\in H^2(M;{\bf Z})|e\in A_{M,\omega} \hbox{ for all $\omega\in \Omega_{M,K}$}\}.$$

Built on results in [T3], [LL1], we are able to compute the
symplectic Seiberg-Witten
invariants on those 4-manifolds and use them to probe a large part of
$A_{M,K}$.
To state the result, we need to introduce yet another concept,
the forward cone associated to a symplectic canonical
class. On a closed oriented  $4-$manifold with
$b^+=1$, the set of classes of positive square ${\cal P}$
fall into two connected
components.
Given an orientation-compatible symplectic form $\omega$,
we will call the  component containing $[\omega]$
the forward cone ${\cal FP}$ associated to $\omega$.
Similarly, given a symplectic canonical class $K$,
we will call the component containing ${\cal C}_K$ the
forward cone associated to $K$ and denote it by
${\cal FP}(K)$.

For a minimal  4-manifold
with $b^+=1$ and a
symplectic canonical class $K$,
the result we need about $A_{M,K}$ is  that (see Proposition
4.2) it contains
large multiples of any class in the forward cone ${\cal FP}(K)$.

Comparing with a result of Donaldson (see [D]), which states that the
Poincar\'e dual of the class of a
sufficiently high multiple of an integral symplectic
form can be represented by symplectic submanifolds,
it is natural to expect that
any class in the forward cone is
represented by symplectic forms.
Using the inflation process of Lalonde and McDuff (see [La], [LM], [Mc3]
 and [Bi2]), we can show that this
is indeed the case.

\noindent{\bf Theorem 2}. Let $M$ be a minimal closed, oriented $4-$manifold
with
$b^+=1$ and $K$ be a  symplectic canonical class.  Then
$${\cal C}_{M,K}={\cal FP}(K).$$
Consequently, any real cohomology class of positive square
is represented by an orientation-compatible symplectic form.

For non-minimal 4-manifolds, it is no longer true
that any real cohomology class of positive square
is represented by an orientation-compatible symplectic form,
due to the presence of the set ${\cal E}$. In fact  what is relavant
for the $K-$symplectic cone is
the subset of $K-$exceptional spheres
$${\cal E}_{M,K}= \{E\in {\cal E}_M|E\cdot K=-1\}.$$
An important fact from [LL2] is that
${\cal E}_{M,K}$ lies in $A_{M,K}$.
From this, together with the blow up formula for
Gromov-Taubes invariants in [LL4],  we can
obtain

\noindent{\bf Theorem 3}.
Let $M$ be a minimal closed, oriented $4-$manifold
with
$b^+=1$ and $K$ be a symplectic canonical class.  Then
$${\cal C}_{M,K}=\{e\in {\cal FP}(K)|e\cdot E>0 \hbox{
 for any $E\in {\cal E}_{M,K}$}\}.$$

Moreover,
 for the symplectic cone
${\cal C}_M$, we have

\noindent{\bf Theorem 4}. Let $M$ be a closed, oriented $4-$manifold with
$b^+=1$ and ${\cal C}_M$ non-empty.
Then
$$C_M=\{e\in {\cal P}|0<|e\cdot E|\hbox { for all $E\in {\cal E}$ }\}.$$

Let us remark that, on a complex surface, our $K-$symplectic cone is
similar to (often larger than) the K\"ahler cone.
Apriori the K\"ahler cone is convex because the sum of two K\"ahler forms
is still a K\"ahler form, while our $K-$symplectic cone may  not have this property
because the sum of two symplectic forms may fail to be a symplectic form.
However, for 4-manifolds with $b^+=1$, it is easy to see from Theorem 3 that
our K-symplectic cone turns out to be  convex.

As we mentioned, the knowledge of the set $A_{M,K}$ is crucial in the
understanding of the moduli space of symplectic forms.
We gradually realize that the set itself is an important invariant of $M$.
In the last part of the paper, we devote ourselves
to studying it.

For the same class of 4-manifolds
above, we are  able to   determine which multiples of $K$
are in $A_{M, K}$. As a pleasant byproduct, we obtain an analogue
of Casteluovo's criterion of rational surfaces.

\noindent{\bf Corollary 2}. Let $M$ be a minimal symplectic 4-manifold
with symplectic canonical class $K$. If $b_1=0$ and $2K\ne 0$,
and $2K$  not in $A_{M,K}$, then $M$ is  $CP^2$ or
$S^2\times S^2$.

We also define the rational $K-$surface cone ${\cal S}^{\bf Q}_{M, K}$ to be
the convex subset of $H^2(M;{\bf Q})$ generated by elements in
$A_{M,K}$. The rational $K-$surface cone is a less refined object
but easier to study.
Our main result about the $K-$surface cone is

\noindent{\bf Theorem 5}. Let $M$ be a  closed, oriented $4-$manifold
with $b^+=1$ and symplectic canonical class $K$.
 Then
$${\cal FP}^{\bf Q}(K)+\sum_{E_i\in {\cal E}_K}{\bf Q}^+ E_i
\subset {\cal S}^{\bf Q}_{M,K}\subset \overline{\cal FP}^{\bf Q}(K)+\sum_{E_i\in {\cal E}_K}{\bf Q}^+ E_i,$$
where ${\cal FP}^{\bf Q}(K)$ and $\overline{\cal FP}^{\bf Q}(K)$
 are the sets of rational classes in ${\cal FP}^{\bf Q}(K)$ and $\overline{\cal FP}^{\bf Q}(K)$ respectively.

On a complex surface, by the Nakai-Moishezon criterion,
a rational type $(1,1)$ cohomology class with positive square
is represented by a K\"ahler form
if and only if it is positive on any irreducible
holomorphic curve.
We notice that, for the 4-manifolds in Theorems 2, 3 and 5,
these theorems imply that a rational cohomology class
with positive square is represented by a symplectic form
if and only if it is positive on any class in
$A_{M,K}$, thus provide a sort of symplectic analogue.

 Actually, if we define the rational $K-$symplectic cone ${\cal C}^{\bf Q}_K$ to be
the set of rational classes in the $K-$symplectic cone, then immediate from the definitions, then inside
$H^2(M;{\bf Q})$,
the rational $K-$surface cone is contained in the
dual of the rational $K-$symplectic cone and vice versa.
Moreover,
Theorems 2, 3 and 5 suggest that there is the following
strong duality,

\noindent{\bf Duality Conjecture}. Let $M$ be a closed, oriented 4-manifold
with $b^+=1$.
Suppose $K$ is a symplectic canonical class, then
 the rational $K-$surface cone and the rational
$K-$symplectic cone are dual to each other.

For a minimal 4-manifold, this conjecture simply asks whether
the rational $K-$surface cone is the closure of the forward cone.
We are able to confirm it for several classes of minimal
4-manifolds (see Proposition 6.5).

After the first draft of this paper appeared, we were informed by P. Biran
 that  weaker versions of our
Theorems 2 and 3 appeared in [Bi2] and [Bi3] (see Theorem 3.2 in page 297
in [Bi3] and Theorem 3.2 B in page 9-10 in [Bi2]). The set `Class C'
that appears in his theorems is a rather general notion. In fact,
following from the results in our paper, it is precisely the set of
4-manifolds $M$ with $b^+=1$ and ${\cal C}_M$ non-empty.

The organization of this paper is as follows.
In \S2, we review the Seiberg-Witten theory and the Gromov-Taubes theory
on symplectic 4-manifolds, in particular those with
$b^+=1$. In \S3,
we present the proof of Theorem 1.
 The $K-$symplectic cone and the symplectic cone are studied in
section 4.  In \S4.1, we deal with minimal 4-manifolds and
prove Theorems 2. In \S4.3, we deal with non-minimal 4-manifolds
and prove Theorems 3-4. We would like to mention that some of the results in
\S4 appeared in [Mc3] in slightly weaker form
(see the remark after Proposition 4.3).
In \S5, we prove Corollary 2.  In \S6, we study the Duality conjecture.

We would like to thank V. Kharlamov, J. Koll\'ar, R. Lee, I. Smith and G. Tian for their
interest in this work. We would also like to thank P. Biran for
some very useful remarks.
This research is partially supported by NSF.

\noindent{\bf Conventions}. In the rest of the paper,
we will make the following simplifications.
An integral cohomology
class is identified with its Poincar\'e dual, and
a complex line bundle is identified with its first Chern class.
A symplectic 4-manifold refers to a pair consisting of a closed oriented
4-manifold $M$ with an orientation-compatible symplectic
form $\omega$.
All the symplectic forms are orientation-compatible
and all the diffeomorphisms are orientation preserving.
All the surfaces are embedded unless specified otherwise.
 We will often drop $M$ from the notations like
$\Omega_M$, $\Omega_{M,K}$, ${\cal C}_{M,K}$, $A_{M, \omega}$, $A_{M,K}$ and
${\cal S}_{M, K}$ where there is no confusion.

\bigskip

\noindent{\bf \S2. The Seiberg-Witten invariants
and the Gromov-Taubes invariants of symplectic $4-$manifolds}
\medskip

In this section we first review the Seiberg-Witten theory
on symplectic $4-$manifolds, in particular those
with $b^+=1$.
Then we  review the Gromov-Taubes theory of
counting embedded symplectic surfaces in symplectic $4-$manifolds.
Finally we review the equivalence between the
two theories on symplectic 4-manifolds with $b^+=1$.

\medskip
\noindent{\bf \S2.1. The Seiberg-Witten invariants of
symplectic $4-$manifolds}
\medskip

In this subsection we review the Seiberg-Witen invariants.
For more details, see e.g. [M] and [S],
and  for the Seiberg-Witten invariants on symplectic $4-$manifolds, see
[T1].

The Seiberg-Witten invariant $SW$ is defined on the
set of Spin$^c$ structures ${\cal SP}$. Associated to
each Spin$^c$ structure ${\cal L}$ is a rank 2 complex vector bundle,
whose determinant line  bundle $c_1({\cal L})$
 is called the determinant bundle
of ${\cal L}$. A useful fact to keep in mind is that ${\cal SP}$
is an affine space modelled on $H^2(M;{\bf Z})$. So, when fixing
${\cal L}$,
 any nonzero class $e\in H^2(M;{\bf Z})$ gives rise to
a Spin$^c$ structure, denoted by ${\cal L}\otimes e$.
And any Spin$^c$ structure is of this form.
The determinant line bundle of ${\cal L}\otimes e$
is related to that of ${\cal L}$ by
$c_1({\cal L}\otimes e)=c_1({\cal L})+2e.$

Fix a riemannian metric $g$  and a real self-dual 2-form ${\mu}$ on $M$.
The Seiberg-Witten equations can be written down for
any Spin$^c$ structure ${\cal L}$ (we refer the readers to [M] for
the actual equations). Denote
${\cal M}(M, {\cal L})$ as the moduli space of the Seiberg-Witten
equations.
For generic pairs $(g, \mu)$, ${\cal M}(M, {\cal L})$
has the nice property that it is a compact manifold of real dimension
$$2d({\cal L})={1\over 4}(2\chi(M)+3\sigma(M))-{1\over 4}
(c_1({\cal L})\cdot c_1({\cal L})).$$
Here $\chi(M)$ and $\sigma(M)$ are the Euler characteristic and
the signature of $M$ respectively. Furthermore, an orientation of the real line $det^+=
H^0(M;{\bf R})\otimes \Lambda ^{b_1}H^1(M;{\bf R})\otimes \Lambda ^{b^{+}}
H^{+}(M;{\bf R})$
naturally orients ${\cal M}(M, {\cal L})$. Here
$b_1$ is the first Betti number of $M$ and  $H^+(M;{\bf R})$ is
a maximal positive subspace of $H^2(M;{\bf R})$.
In addition, there is a naturally defined
circle bundle on $M\times {\cal M}(M, {\cal L})$, whose Euler class induces
a  map
$\phi$ from $H_{*}(M;{\bf Z})$ to $H^{2-*}({\cal M}(M, {\cal L});{\bf Z})$.

Now we give the definition of $SW$:

\noindent {\bf Definition 2.1}. Fix an orientation for the line
$det^+$ and a Spin$^c$ structure ${\cal L}$. Fix also  a generic
pair $(g, \mu)$ with the salient properties above. Let
$\gamma_1\wedge\cdots\wedge\gamma_p$ be a decomposable element in
$\Lambda^p(H_1(M;{\bf Z})/Torsion)$. Then the value of $SW (M,
{\cal L})\in \Lambda^*H^1(M;{\bf Z})$ is defined as follows:

$$ SW(M,{\cal L}; \gamma_1\wedge\cdots \wedge \gamma_p)=
\int_{{\cal M} (M,{\cal L})}
\phi(\gamma_1)\wedge\cdots\wedge\phi(\gamma_p)\wedge
\phi(pt)^{(2d-p)/2}, \eqno (2.1) $$
where
$pt$ is the class of a point in $H_0(M;{\bf Z})$.

When $b^+>1$, the value of $SW$ is independent of the choice of
the generic pair $(g, \mu)$. Hence $SW$ can be viewed as a map from
the set of Spin$^c$ structures ${\cal SP}$ to $\Lambda^*H^1(M;{\bf Z})$.
Let $SW^i(M, {\cal L})$ denote the part of $SW(M,{\cal L})$ in
 $\Lambda^i H^1(M;{\bf Z})$.
Suppose $\eta$ is the positive generator of
$\Lambda^0 H^1(M;{\bf Z})\equiv {\bf Z}$.
Then the integer $SW^0(M, {\cal L}; \eta)$ will be
simply denoted by $SW^0(M,{\cal L})$.

When $b^+=1$, $SW$ also depends  on the choice of a chamber.
On a 4-manifold with $b_2^+=1$, the set of real
second cohomology classes with positive
square ${\cal P}$ is
 a cone with two connected components. Pick one of them and
call it the forward cone. Given a metric $g$, there is a unique
self-dual harmonic 2-form $\omega_g$ for $g$ in the forward cone with
$\omega_g^2=1$. For a pair $(g, \mu)$ and a Spin$^c$
structure ${\cal L}$, define the discriminant $\triangle_{\cal L}(g, \mu)=
\int (c_1({\cal L})-\mu)\omega_g$.
 The set of pairs $(g, \mu)$ with positive and negative
discriminant are called the positive and negative ${\cal L}$
chamber respectively. The map $SW(M,{\cal L})$ is constant on any
${\cal L}$ chamber. So in the case $b^+=1$, given a choice of the
the forward cone and a Spin$^c$ structure, we can define $SW_{+}(M,{\cal L})$
and $SW_{-}(M, {\cal L})$ for the positive and negative ${\cal L}$
chambers respectively. $SW^i_{\pm}(M,{\cal L})$
are similarly defined.

 From now on, when there is no confusion, we will often drop $M$ from
  $SW^i(M, {\cal L}_{K^{-1}}\otimes e)$ or $SW^i_{\pm}(M, {\cal L}_{K^{-1}}\otimes e)$.

Recall that  a symplectic 4-manifold is a closed oriented
4-manifold
with an orientation-compatible symplectic form $\omega$.
When $M$ is a symplectic 4-manifold, we have the
following facts.

\noindent 1. As mentioned in \S1,
the symplectic form $\omega$ (actually, the deformation class of $\omega$)
determines a unique homotopy class of
$\omega-$compatible almost complex structures, and hence a canonical
line bundle $K$.
These almost complex structures induce
a Spin$^c$ structure ${\cal L}_{K^{-1}}$ with   $K^{-1}$ as
its determinant line bundle.

\noindent 2. Since $K\cdot K=2\chi(M)+2\sigma(M)$,
  the formal dimension of the Seiberg-Witten moduli space
 of ${\cal L}_{K^{-1}}\otimes e$ is then give by
$$2d=e\cdot e-K\cdot e.\eqno (2.2)$$

\noindent 3. There is a natural orientation of the line
$det^+$.

\noindent 4.  When $b^+=1$,
as mentioned in \S1, $\omega$ determines the choice of the forward cone.
 With respect to this choice, the negative chamber is called the
 symplectic chamber.

Hence we will use
the symplectic structure to orient the Seiberg-Witten moduli spaces. In the case
$b^+>1$, the SW invariant so defined will be denoted by $SW_{\omega}$;
and, in the case $b^+=1$, use the symplectic structure to define $SW_{\omega,+}$ and $SW_{\omega, -}$ as well.

 A fundamental result of Taubes is

\noindent {\bf Theorem 2.2} (see [T1]). Let $M$ be a
symplectic
$4-$manifold with symplectic canonical class $K$.
Then $SW_{\omega}^0({\cal L}_{K^{-1}})=1$ when $b^+>1$, and
$SW^0_{\omega,-}({\cal L}_{K^{-1}})= 1$ when $b^+=1$.

An involution on the set of Spin$^c$ structure induces the following
symmetry of the Seiberg-Witten invariants,

\noindent{\bf Symmetry Lemma 2.3} ([W]). Let $M$ be a symplectic
$4-$manifold with $b^+=1$ and symplectic canonical class $K$.
Then, for any integral class $e$,
$$SW^i_{\omega,+}({{\cal L}_{K^{-1}}\otimes e})=(-1)^{(1-b_1-i+b^+)/2}
SW^i_{\omega,-}({\cal L}_{K^{-1}}\otimes (K-e)).$$

Notice that, by the symmetry lemma and Theorem 2.2, we have  $SW^0_+({\cal L}_{K^{-1}}\otimes
K)=\pm 1$.

The following inequalities constrain the range of symplectic forms.

\noindent{\bf Theorem 2.4}. Let $M$ be a smooth,
closed oriented $4-$manifold.

\noindent 1 ([T2]). Suppose $M$ has $b^+>1$ and $e$ is represented by a symplectic form $\omega$ with
$K$ as its symplectic canonical class. If  $L$ is an integral  class such that
$SW_{\omega}({\cal L}_{K^{-1}}\otimes L)$ is nontrivial,
then $K\cdot e\geq L\cdot e$ with equality only if $K=L$.

\noindent 2 ([LL2]). Suppose $M$ has $b^+=1$. Suppose $e_1$ and $e_2$ are
 represented
by symplectic forms  with $K_1$ and $K_2$ as their symplectic canonical classes respectively, and satisfies $e_1\cdot e_2>0$.
Then $K_1\cdot e_1\geq K_2\cdot e_1$,  and equality holds
only if $K_1-K_2$ is a
torsion class.

\medskip
\noindent{\bf \S2.2. Gromov-Taubes invariants of symplectic $4-$manifolds}
\medskip
Let $M,\omega$ be a closed symplectic $4-$manifold with symplectic canonical
class $K$.
Like the Seiberg-Witten invariants, the Gromov-Taubes invariants,
introduced  by Taubes in [T3], are defined
for any class $e\in H^2(M;{\bf Z})$ and takes values in
 $\oplus \Lambda ^*H^1(M;{\bf Z})$.

When  $e$ is
a nonzero integral class, introduce
 the Gromov-Taubes dimension:
$$d(e)=e\cdot e-K\cdot e.$$
$d(e)$  is the expected maximal complex dimension of the
components of pseudo-holomorphic curves (the domain of the curves
can be any Riemann surfaces with arbitrary number of connected
components) representing  $e$.

If $d(e)\geq 0$, fix an integer $p\in\{0,1,\cdots,d\}$ and then fix
$\gamma_1\wedge \cdots \wedge \gamma_{2p}
\in \Lambda^{2p}H_1(M;{\bf Z})$.
As in the case of Seiberg-Witten invariants,  it suffices to
define $Gr_{\omega}(M, e)(\gamma_1\wedge \cdots \wedge \gamma_{2p})$.
Fixing a
compatible
almost complex structure $J$,
let
${\cal H}(e,J,Z, \Gamma)$
 be the set
of
 $J-$holomorphic curves representing
  $e$, and passing through
a set of
$d-p$
points  $Z$, and intersecting $\Gamma$, a set of $2p$  disjoint
circles representing
$\gamma_1, \cdots, \gamma_{2p}$.

For generic $J, Z$  and $\Gamma$,
${\cal H}(e,J,Z, \Gamma)$ has the following properties:

\noindent 1. ${\cal H}(e,J,Z, \Gamma)$
is a finite set.

\noindent 2. Let $h\in {\cal H}(J, Z, \Gamma)$. Let $C_1, \cdots,C_k$
be the irreducible components of $h$ which represent
 the classes $e_1,\cdots , e_k$
and have multiplicities $m_1,\cdots, m_k$. Then $C_1,\cdots, C_k$
are embedded and disjoint.

\noindent 3. $m_j=1$ unless $C_j$ is a torus and $e_j^2=0$.

\noindent 4. $e_j^2\geq -1$, and $e_j^2=-1$ only if $C_k$ is a sphere.

\noindent 5.  To each $h\in {\cal H}(e, J, Z, \Gamma)$, an integer $q(h)$  can be assigned
in a delicate way.

\noindent{\bf Definition 2.5}.
When $e$ is the zero class, $Gr_{\omega}(M,0)$ is simply defined to
be 1. When $e$ is not zero, fix
$\gamma_1\wedge \cdots \wedge \gamma_{2p}
\in \Lambda^{2p}H_1(M;{\bf Z})$.
Choose
generic $(J, Z, \Gamma)$ with the salient properties above.
$Gr_{\omega}(M, e)(\gamma_1\wedge \cdots \wedge \gamma_{2p})$
is then defined to be $\sum q(h)$.

When there is no confusion, we will freely
drop $M$ from  $Gr_{\omega}(M, e)$.

The Gromov-Taubes invariants
are independent of the choice of the generic $(J, Z, \Gamma)$.
In fact they
 only depend on the deformation class of the symplectic form
and it is natural with respect to diffeomorphisms.
A version of Gromov-Taubes invariants was introduced by  Ruan
in [R3]. It was also shown in [IP] that the Gromov-Taubes invariants
can be constructed from the Ruan-Tian invariants in [RT].

Suppose $e$ is a class with non-trivial Gromov-Taubes invariant.
Then for generic $J, Z, \Gamma)$, there
exists a  $J-$holomorphic curve  $h$ representing
$e$ satisfying the properties 2-4 above.
Since any embedded $J-$holomorphic curve is an embedded symplectic surface,
it is clear that, if $m_j=1$ for all $j$, then
$e$ is represented by an embedded symplectic surface.
If $m_j>1$ for some $j$, then $C_j$
is a torus with square zero.
In this case, we have
the following simple fact.

\noindent{\bf Lemma 2.6}. Let $M$ be a symplectic $4-$manifold and
$e$ be an integral  class with square zero. If
$e$ can be represented by
connected embedded symplectic surfaces, then  a positive multiple of $e$
can also be represented by embedded symplectic surfaces.

\noindent{\it Proof}. Let $C$ be a connected embedded symplectic surface
 representing $e$. By the symplectic neighborhood theorem, the tubular neighborhood
of $C$ is symplectically a product $C\times D^2$, where
$D^2$ is a 2-disk. For any positive
integer $n$, pick $n$
points in $D^2$ and we get $n$ disjoint embedded
symplectic surfaces whose disjoint union
represents $ne$. The proof is complete.

Thus by Properties 1$-$4 and Lemma 2.6 we have

\noindent{\bf Theorem 2.7} ([T4]). Let $M$ be a symplectic
$4-$manifold. If $e\in H^2(M;{\bf Z})$ and
$Gr_{\omega}(e)$
is non-trivial, then
 $e$ can be represented by an embedded symplectic surface
whose only components with negative square are
spheres with square $-1$.

Given a symplectic sphere $\Sigma$ with square -1 in $M$, one
can replace a neighborhood of $\Sigma$ by a symplectic
ball to obtain a new  symplectic 4-manifold. This process
is called  (symplectic) blowing down. The reverse process is
called (symplectic) blowing up (at a point in the symplectic
ball).

Let ${\cal E}_{M, \omega}$ be the set of classes which
are represented by
symplectic $-1$ spheres. $M$ is called symplectically minimal if
${\cal E}_{M, \omega}$ is empty.

Since every symplectic 4-manifold $M$ can be obtained by blowing
up a minimal symplectic 4-manifold $N$ at a number of points (see
[Mc4]), the following blow up formula is useful to  compute the
Gromov-Taubes invariants of non-minimal symplectic 4-manifolds.

\noindent{\bf Theorem 2.8} (see [LL4]).
 Let $M,\omega$ be a symplectic 4-manifold with $b^+=1$.
Suppose $E\in {\cal E}_{M,\omega}$  is represented
by the symplectic sphere $\Sigma$.
Let $N,\beta$ be the symplectic manifold obtained by blowing
down $\Sigma$. If $v\in H^2(N;{\bf Z})$
and $u=v-lE$ for some integer $l$.
Then, under the canonical identification between
$\Lambda ^{*}H^1(N;{\bf Z})$ and $\Lambda ^{*}H^1(M;{\bf Z})$,
$$Gr_{\omega}(M, u)=Pr(2d(u))Gr_{\beta}(N, v).$$
Here  $pr(2d(u))$ is the projection from
$ \Lambda^*H^1(M;{\bf Z})$ to $ \oplus_{i=0}^{2d(u)}
\Lambda^i H^1(N;{\bf Z})$.

\medskip
\noindent{\bf \S2.3. Equivalence between SW and Gr in the case $b^+=1$}
\medskip

In this subsection, $M$ is a symplectic $4-$manifold with $b^+=1$.
Another fundamental and deep result of Taubes ([T3]) is

\noindent{\bf Theorem 2.9}.
Let $M,\omega$ be a symplectic 4-manifold with $b^+=1$.

\noindent 1. If $SW_{\omega,-}({\cal L}_{K^{-1}}\otimes e)$ is non-trivial, then
there are positive integers $n_i$ and classes $e_i$ such that
$e$ can be written as $e=\sum_i n_i e_i$ and
$e_i$ is represented by an $\omega-$symplectic surface.
In particular, if $SW_{\omega,-}({\cal L}_{K^{-1}}\otimes e)$ is nonzero, then $e\cdot \omega\geq 0$
and $e\cdot \omega=0$ only if $e=0$.

\noindent 2. If $M$ is minimal,
then $SW_{\omega,-}=Gr_{\omega}$.
 For a non-minimal symplectic
4-manifold, the above conclusion holds with
the additional assumption that

$$ e\cdot E\geq -1 \hbox{   for each  }E\in {\cal E}_{\omega}\eqno (2.2).$$

For classes violating the condition (2.2), $SW_{\omega,-}$ and
$Gr_{\omega}$ are different. To extend the
equivalence to those classes
McDuff introduces a map $Gr_{\omega}'$ (see [Mc2]). It is a variation of the map $Gr$ which
takes into account
multiply covered $-1$ spheres,
and coincides with $Gr_{\omega}$ for all classes satisfying (2.2). It
was shown in [LL4] that indeed $Gr_{\omega}'=SW_{\omega,-}$ for all classes not satisfying (2.2).

\bigskip
\noindent{\bf \S3. Uniqueness of symplectic canonical class}
\medskip

In this section we will present the proof of
Theorem 1 and Corollary 1, which rely on
the computation of the Seiberg-Witten invariants.

We start with citing the following simple lemma concerning the intersection
pairing
of a $4-$manifold with $b^+=1$ is very useful (see [LL2]).

\noindent{\bf Lemma 3.1} (Light cone lemma).
Suppose $M$ is a  manifold
with $b^+=1$. Let $A$ and $B$ be 2 classes in $H^2(M;{\bf R})$ with $A^2\geq 0$.

\noindent 1. If $B\cdot A=0$, then $B^2\leq 0$. And $B^2=0$ iff $A^2=0$ and
$B=rA$ up to torsion.

\noindent 2. If $A$ and $B$ are all in the closure of the forward cone, then
$A\cdot B\geq 0$.

Now we turn to the computation of the SW invariants.

\noindent{\bf Lemma 3.2}. Let $M,\omega$ be a  minimal symplectic manifold
with
$b^+=1$ and the symplectic canonical class $K$. Let $e$ be a class in $H^2(M;{\bf Z})$.

\noindent I. Suppose $K^2\geq 0$ and  $K\cdot \omega\geq 0$.

\noindent If  $({K^{-1}}+2 e)\cdot \omega\geq 0$ and  $SW_{\omega,+}({\cal L}_{K^{-1}}\otimes e)$  is non-trivial, then

\noindent I.1. $e=K$ if $K^2>0$ or $K$ is a torsion class.

\noindent I.2. $e=rK$ for some rational number $r\in [1/2, 1]$ or
$2e-K$ is a torsion class  if $K^2=0$.

\noindent II. Similarly, if  $({K^{-1}}+2 e)\cdot \omega \leq 0$ and
$SW_{\omega,-}({\cal L}_{K^{-1}}\otimes e)$ is nontrivial, then

\noindent II.1. $e$ is the zero class if $K^2>0$ or $K$ is a torsion class.

\noindent II.2. $e=rK$ for some rational number $r\in [0, 1/2]$
or $2e-K$ is a torsion class if $K^2=0$.

\noindent{\it Proof}.
Since
$K^2\geq 0$ and $K\cdot \omega\geq 0$, $K$ is in the closure of the
forward cone.
%In either cases, $e\cdot e-K\cdot e$, which
%is the Seiberg-Witten dimension
%of ${\cal L}_{K^{-1}}\otimes e$, is non-negative.

\noindent{\it Case I}:
  $({\cal L}_{K^{-1}}+2e)\cdot \omega\geq 0$ and  $SW_{\omega,+}({\cal L}_{K^{-1}}\otimes e)$  is non-trivial.

By the Symmetry lemma 2.3,
$$SW^i_{\omega,+}({\cal L}_{K^{-1}}\otimes e)=\pm
SW^i_{\omega,-}({\cal L}_{K^{-1}}\otimes 2(K-e)).$$
Thus
$SW_{\omega,-}({\cal L}_{K^{-1}}\otimes 2(K-e))$ is non-trivial.
By Theorem 2.9, if $K\ne e$, then
 $(K-e)\cdot \omega >0$ and $Gr_{\omega}(K-e)$ is non-trivial. Furthermore,
since $M$ is minimal, by Theorem 2.7,  $K-e$ is represented
by an embedded symplectic surface whose  components all have
non-negative square.  So  $(K-e)^2\geq 0$ and
$K-e$ is in the forward cone.

Let us first deal with the case that $K$ is a torsion class.
Since $K\cdot \omega=0$, we see that  $e\cdot \omega\geq 0$. And
if $e\ne K$, $e\cdot \omega=(e-K)\cdot \omega<0$.
This is a contradiction. Therefore, in this case, $e$ must be equal to $K$.

Now let us treat the case that $K$ is not a torsion class.
 Since
$K-e$ is in the forward cone,
$$(K-e)\cdot K\geq 0 \eqno (3.1) .$$

Suppose $({K^{-1}}-2e)\cdot \omega=0$, then $K-2e$ is a torsion class.
The Seiberg-Witten dimension of ${\cal L}_{K^{-1}}\otimes e$,
which should be positive,  is then given by
$e\cdot e-K\cdot e=-K^2/2$.  Since $K^2$ is assumed to be non-negative,
 we must have
 $K^2=0$.

Now assume  that $(K^{-1}+2e)\cdot \omega$ is strictly positive.
So there is a $c\geq 1$ such that $(K-c(K^{-1}+2e))\cdot \omega=0$
and hence
$(K-c(K^{-1}+2e))^2\leq 0$, with equality only if $K-c(K^{-1}+2e)=0$.
 However,
$$\eqalign{(K-c(K^{-1}+2e))^2&=((1+c)K-2ce)^2\cr
&=(1-c)^2K^2    +4c(K^2-e\cdot K)+  4c^2(e^2-e\cdot K).\cr}$$
All three terms are non-negative, the first by assumption, the second by
(3.1), the third being the Seiberg-Witten dimension of  ${\cal L}_{
K^{-1}}\otimes e$.

When $K^2>0$,
equality holds if and only if $c=1$, and thus $(K-e)\cdot \omega=0$.

When $K^2=0$, equality holds if and only if  $e^2=e\cdot K=0$. Since $K$
is not a torsion class, by the light cone lemma 3.1,
$e=rK$. $r$ is no bigger than one since $(K-e)\cdot \omega\geq 0$ and
$r$ is no less than $1/2$ since $(K^{-1}+2e)\cdot \omega\geq 0$.
We have thus finished the proof of case I.

\noindent{\it Case II}:
  $(K^{-1}+2e)\cdot \omega\leq 0$ and  $SW_{\omega,-}({\cal L}_{K^{-1}}\otimes e)$  is
non-trivial.

Then $(K^{-1}+2(K-e))\cdot \omega\geq 0$ and $SW_{\omega,+}({\cal L}_{K^{-1}}\otimes
(K-e))$ is
non-trivial. Thus we can similarly determine $K-e$ and hence $e$ itself.
The proof of the lemma is complete.

We will need  the wall crossing formulas in [LL1], [LL4]
(see also [OO] and [OT]).

\noindent{\bf Lemma 3.3}. Let $M,\omega$ be a symplectic
4-manifold
with $b^+=1$ and symplectic canonical class
$K$. Let $e$ be a class in $H^2(M;{\bf Z})$.
When $b^+=1$,  the choice of the forward cone is just an orientation of
the
line $H^{+}=\Lambda ^{b^{+}}H^{+}$, so it induces an orientation
of $\Lambda^{b_1}H^1$.

%In particular $SW_{\omega,-}(K^{-1}\otimes e)-SW_{\omega,+}(K^{-1}\otimes e)$ is nonzero if

\noindent 1. Suppose $d(e)\geq b_1$. Let $\gamma_1\cdots, \gamma_{b_1}$ be a basis
of $H_1(M;{\bf Z})/Torsion$ such that $\gamma_1\wedge\cdots\wedge\gamma_p$ is the dual
orientation of the symplectic orientation on $\Lambda^{b_1}H^1(M;{\bf Z})$.
Then
$$\eqalign{&SW^{b_1}_{\omega,-}({\cal L}_{K^{-1}}\otimes e;
\gamma_1\wedge\cdots\wedge\gamma_p)-SW^{b_1}_{\omega,+}({\cal L}_{K^{-1}}\otimes e;
\gamma_1\wedge\cdots\wedge\gamma_p)=1.\cr}$$

\noindent 2. Suppose $d(e)\geq 0$ and $\Lambda^2H^1(M;{\bf Z})$ is
of rank 1.
Let $\gamma$ be the generator of $\Lambda^2 H^1(M;{\bf Z})\subset H^2(M;{\bf Z})$ such that
$\omega \cdot \gamma >0$, then

$$SW^0_{\omega,-}({\cal L}_{K^{-1}}\otimes e)-SW^0_{\omega,+}({\cal L}_{K^{-1}}\otimes e)={(K^{-1}+2e)\cdot \gamma \over 2}.$$

%\noindent{\it Proof}. We only have to fix the sign.
%This can be accomplished by checking special manifolds.

\noindent{\bf Lemma 3.4}. Let $M$ be a
rational or irrational ruled symplectic $4-$manifold with
symplectic canonical class $K$. For any class $e\in H^2(M;{\bf Z})$,
$SW^0_{\omega,-}({\cal L}_{K^{-1}}\otimes e)$ is non-trivial if $(K^{-1}+2e)$ is in the
forward cone and $d(e)\geq 0$.

\noindent{\it Proof}. Choose a metric $g$ of positive scalar curvature
on $M$. With such a metric and the zero self-dual 2-form, the
Seiberg-Witten moduli space of each Spin$^c$ structure
is empty ([W]). Since $K^{-1}+2e$ is assumed to be in the forward cone,
$(K^{-1}+2e)\cdot \omega_g>0$,
where $\omega_g$ is the unique self-dual harmonic
2-form in the forward cone. So the pair $(g,0)$ is in the positive
chamber
for the Spin$^c$ structure ${\cal L}_{K^{-1}}\otimes e$.
 We thus find that
$SW_{\omega,+}({\cal L}_{K^{-1}}\otimes e)=0$.
For a rational 4-manifold, $b_1=0$, the conclusion then follows from part
one of
Lemma 3.3. For an irrational ruled 4-manifold, $b_1\geq 2$ and $\gamma$ is
nonzero.  Since $(K^{-1}+ 2e)$ is in the forward cone and $\gamma$ has
square zero, $(K^{-1}+2e)\cdot \gamma$ is nonzero. Hence the conclusion
in this case
follows from part two of Lemma 3.3.

Now we need to review the notion of minimal reduction.
Recall we  defined ${\cal E}_M$ in \S1, and  $M$ is said to be
(smoothly)
minimal if ${\cal E}_M$ is empty. Obviously,
if $M$ is a symplectic 4-manifold, ${\cal E}_{M,\omega}$ is a subset of ${\cal E}_M$,
and so $M$ is symplectically minimal if it is
(smoothly) minimal.   In fact it follows from
the results in [L1] that the reverse is also true.

 Any 4-manifold $M$ can be
decomposed as a connected sum of a minimal manifold $N$  with some
number of $\overline {CP}^2$. Such a decomposition is called
a (smooth) minimal reduction of $M$, and $N$ is
a minimal model
of $M$. $M$ is
said to be
 rational if one of its minimal models is $CP^2$ or $S^2\times S^2$;
and irrational ruled if one of its minimal models is an
 $S^2-$bundle over a Riemann surface of positive genus.
When $M$ has non-empty symplectic cone and
is not rational nor irrational ruled,
$M$ has a unique minimal reduction (see [L1] and also [Mc3]).

The only minimal rational manifolds are $CP^2$ and $S^2\times S^2$.
 And  a non-minimal rational manifold  has two kinds of decompositions.
It is either decomposed as $CP^2\# l\overline {CP}^2$ or as
$S^2\times S^2 \# (l-1)\overline {CP}^2$.
We will always use the first decomposition and call it a standard
decomposition.
 The picture for irrational ruled manifolds is
similar.
$S^2-$bundles over a Riemann surface of positive genus are the only
minimal irrational ruled manifolds. Fix the base  surface $\Sigma_g$,
there are two $S^2-$bundles over it, the trivial one
$S^2\times \Sigma_g$ and the unique
non-trivial one $S^2\tilde{\times}\Sigma_g$.
A non-minimal  irrational ruled manifold also has two types of
decompositions.   It is either decomposed as
$S^2\times \Sigma_g\# l\overline {CP}^2$ or as
$ S^2\tilde{\times}\Sigma_g \# l\overline {CP}^2$. We will use the first
decomposition and call it a standard decomposition.

Let $H$ be a generator of $H^2(CP^2;{\bf Z})$ and $F_1,\dots, F_l$ be
the generators of  $H^2$ of the
$\overline{CP}^2$.
Let $U$ and $T$ be  classes in $S^2\times \Sigma_g$ represented by
$\{pt\}\times \Sigma_g$  and  $S^2\times \{pt\}$ respectively.
$H, F_1,\dots, F_l$ are naturally considered as classes in
$H^2(CP^2\# l\overline {CP}^2;{\bf Z}) $ and form a basis.
We will call such a basis  a standard basis. Similarly,
$U, T, F_1, \dots, F_l$ are naturally considered as classes
in  $H^2(S^2\times \Sigma_g\# l\overline {CP}^2; {\bf Z})$ and form a basis.
Such a basis is also called a standard basis.
 On $CP^2\# l\overline {CP}^2$, let  $K_{0}=-3H+\sum_i F_i$;
and on $S^2\times \Sigma_g\# l\overline {CP}^2$, let  $K_0=2U+(2-2g)T+\sum_i F_l$. By the blow up construction (see e.g. [Mc1])
$K_0$ is a symplectic canonical class.

For a  4-manifold $M$ and a choice of symplectic canonical class $K$,
recall the set of $K-$exceptional spheres introduced in \S1 is
${\cal E}_K=\{E\in {\cal E}|K\cdot E=-1\}.$ We will need
the following facts about ${\cal E}_K$.

\noindent{\bf Lemma 3.5}.  Let $M, \omega$ be a symplectic 4-manifold with
$K$ as its symplectic canonical class.

\noindent 1. Let $f$ be a diffeomorphism, then
$ {\cal E}_{f^*K}=
f^* {\cal E}_K.$

\noindent 2.
${\cal E}_K={\cal E}_{\omega},$ i.e. every
class in ${\cal E}_K$ is represented by
an $\omega-$symplectic exceptional sphere.
Moreover, if $E_1, ..., E_p\in {\cal E}_K$
are pairwise orthogonal, then they are represented
by $p$ disjoint $\omega-$symplectic spheres.

\noindent 3. If $E_1$ and $E_2$ are two distinct classes in ${\cal E}_K$,
then $E_1\cdot E_2\geq 0$.

\noindent 4. Suppose
$M$ is not rational nor irrational ruled and
it has a smooth minimal
reduction $N\# l\overline{CP}^2$.
with $F_i$   a generator of
$H^2$ of the $i-$th
$\overline {CP}^2$.
  If we let
$\delta_i=K\cdot F_i$, then $\delta_i=\pm 1$ and
$${\cal E}_K=\{-\delta_1 F_1, ..., -\delta_l F_l\},$$

\noindent{\it Proof}.
$C$ is a $-1$ sphere symplectic with respect to $\omega$ if and only if
$f^{-1}(C)$ is a $-1$ sphere symplectic with respect to $f^*\omega$.
Since  the (cohomology) class represented by $f^{-1}(C)$ is the pull back of
 the
(cohomology) class represented  $C$, and the symplectic canonical class
of $f^*\omega$ is $f^*K$, we have
$f^* {\cal E}_K={\cal E}_{f^*K}$.

 2 was proved in [LL2] and [T3], and 3 follows from 2 and the
fact proved in [Mc1] that, if $E_1$ and $E_2$ are in ${\cal E}_{\omega}$
and distinct,
then $E_1\cdot E_2\geq 0$.
Finally  4 was proved in [L1].

\noindent{\bf Proposition 3.6}. Suppose $M, \omega$ is a
symplectic 4-manifold with symplectic canonical class $K$.
Suppose it has a (smooth) minimal
reduction $N\# l\overline{CP}^2$.
 Let $F_i$  be a generator of
$H^2$ of the $i-$th
$\overline {CP}^2$.
Then
there is a symplectic form $\beta$ on $N$,
such that $(M, \omega)$ is obtained by blowing up
$(N, \beta)$ at $l$ points.
Moreover, there is a diffeomorphism
carrying $K$ to $V\pm F_1+....\pm F_l,$ where
$V$ is a symplectic canonical class of $N$.

\noindent {\it Proof}.  When there are $l$ disjoint
$\omega-$symplectic $-1$ spheres in $M$, we can simultaneously
blow them down  to  obtain a symplectic 4-manifold $N'$. Suppose
these $\omega-$symplectic spheres represent $F_1', ..., F_l'$ and
$V'$ is the symplectic canonical class of $N'$. Then
$K=V'+F_1'+...+F_l'$. And if there is a diffeomorphism $\phi$ such
that $\phi^* F_i'=\pm F_i$ for each $i$, then $N'$ is
diffeomorphic to $N$.

We start with the easier case when $M$ is not rational nor ruled.
By Lemma 3.5,  ${\cal E}_{\omega}=\{\delta_1 F_1,..., \delta_l
F_l\},$ where $\delta_i=\pm 1$. Since $F_i\cdot F_j=0$ for each
pair $i\ne j$, we can find $l$ disjoint $\omega-$symplectic
spheres representing $\delta_1F_1, ..., \delta_l F_l$. Blowing
down these spheres, we obtain a symplectic $4-$manifold $N,
\beta$. If $V$ is the symplectic canonical class of $\beta$, then
$K=V+\delta_1 F_1+...+\delta_l F_l.$

Suppose $M$ is rational or irrational ruled.
The $F_i$ are represented by disjoint smoothly embedded $-1$ spheres.
By the proof of Theorem 1 in [L1], there exists a diffeomorphism $\phi_l$ such that
$F_l'=\phi_l^*F_l$ satisfies $F_l'\cdot K=-1$.

To describe the diffeomorphism, recall that
if $\alpha$ is a class with
square $-1$, then one can define an automorphism of
$H^2$,
the reflection $R(\alpha)$ along $\alpha$, as follows:
$$R(\alpha)\beta=\beta+2(\beta\cdot \alpha)\alpha.$$
And if $\alpha\in {\cal E}$, then the reflection is
realized by a diffeomorphism.

$\phi_l$ is in fact constructed as a composition of
reflections along a series of classes $Y_1, ..., Y_p$
represented by $\omega-$symplectic -1 spheres.
 Moreover, if we go through the proof carefully and use the
fact, mentioned in the proof of
Lemma 3.5,  that two $\omega-$symplectic -1 spheres in two distinct
classes have non-negative intersection,
we find that the classes $Y_i$ have the
following property: if $A$ is represented by an $\omega-$symplectic
-1 sphere and $A\cdot F_l=0$, then $A\cdot Y_i=0$ for each $i$.
Therefore $A$ is invariant under $\phi_l$.

%is represented by an $\omega-$symplectic $-1$ sphere $C_l$.
%We can blow down $C_l$ to obtain a symplectic manifold diffeomorphic to
%(and denoted by) $N_{l-1}$.
The class $\phi_l^*E_{l-1}$
%are naturally in $H^2(N_{l-1};{\bf Z})$
%and since blowing down is a local operation,  they are
is still represented by a
smoothly embedded $-1$ sphere.
By repeating the above process, we find a diffeomorphism
$\phi_{l-1}$ such that $\phi_{l-1}^*\phi_l^*E_{l-1}\cdot K=-1$.
Moreover, since $F_l'=\phi_l^*F_l$ is orthogonal to
$\phi^*E_{l-1}$ and  represented by
an $\omega-$symplectic -1 sphere, $\phi^*_{l-1}F_l'=F_l'$.
Repeating this process $l-2$ more times, we get
$l-2$ more diffeomorphisms $\phi_{l-2},..., \phi_1$
such that
the diffeomorphism
$\phi=\phi_1\circ...\circ\phi_l$ has the property that
$\phi^*F_i\cdot K=-1$ for each $i$.
Therefore we find that  the  symplectic 4-manifold $M,\omega$ can
indeed be blown down to a symplectic $4-$manifold diffeomorphic
to  $N$.

To prove the last statement,
consider  $\omega'=\phi^*\omega$. Its symplectic canonical
class is $\phi^*K$ and $F_i\cdot \phi^*K=-1$.
So we can blow down $l$ disjoint $\omega'-$
symplectic -1 spheres in $M, \omega'$ to obtain a symplectic
$4-$manifold diffeomorphic to $N$. If $V$ is the symplectic canonical
class, then $\phi^*K=V+F_1+...+F_l$, i.e $K$ is carried by $\phi$ to
$V=F_1+...+F_l$. The proof is complete.

We are ready to prove Theorem 1.

\noindent{\bf Theorem 1}. Let $M$ be a smooth, closed oriented manifold
with
$b^+=1$.
There is exactly one equivalence class of symplectic canonical classes.
In fact, if $M$ is minimal,  the symplectic canonical class is unique up to
sign.

\noindent{\it Proof}. We first prove that, when
$M$ is minimal, if $K$ is a symplectic canonical class,
then
the only other symplectic canonical class is $-K$.

\noindent{\it Case 1. $M$ not rational nor irrational ruled}.

\noindent Fix an orientation-compatible symplectic structure
$\omega$ and
let $K$ be its symplectic canonical class.  We have
$K\cdot\omega\geq 0$
and $K^2\geq 0$ by [Liu]. Suppose $\tilde \omega$ is another orientation-compatible
symplectic form and $\tilde
K$
is its symplectic canonical
class. Since the symplectic canonical class of $-\tilde \omega$ is
just $-\tilde K$, we can
assume that $\tilde \omega$ and $\omega$ are
in the same component of the positive cone so they determine the
same forward cone.
Again we have $\tilde K\cdot \omega\geq 0$ and $SW^0_-( {\cal L}_{\tilde
K^{-1}})=1$.
Thus $SW^0_+({\cal L}_{\tilde
K^{-1}}\otimes \tilde K)=\pm 1$ by the symmetry lemma 2.3.

The Spin$^c$ structure
${\cal L}_{\tilde
K^{-1}}\otimes \tilde K$ is of the form
 ${\cal L}_{K^{-1}}\otimes e$ for some class $e$.
Comparing the determinant line bundles, we have
  $$\tilde K^{-1}+2\tilde K=K^{-1}+2e,$$
and so $2e=K+\tilde K$. To prove that $K=\tilde K$, it suffices to show that
$e=K$.

Since $\tilde K$ is in the closure of the forward cone,
we have $(K^{-1}+2e)\cdot \omega \geq 0$.
By the first part of Lemma 3.2, if $K^2>0$ or $K$ is a torsion class, $e$ is equal to $K$. Therefore $K=\tilde K$ in this case.

By the second part of Lemma 3.2, in the remaining case when $K^2=0$ and $K$ is
not a torsion class,
either  $e=r K$ for some rational number $r\leq 1$
or $2e-K$ is a torsion class.
Since $2e=K+\tilde K$, we have either $\tilde K=(2r-1)K$ for some rational
number $r \in [1/2, 1]$, or $2\tilde K-K$ is a torsion class.
In both cases, it is easy to see that
$\tilde K^2=0$ and $\tilde K$ is not a torsion class.  Now if we
start with $\tilde \omega$ and $\tilde K$ and repeat the argument above,
we conclude that
either $K=(2\tilde r-1)\tilde K$ for some rational number $\tilde r \in
[1/2, 1]$, or $2K-\tilde K$ is a torsion class.
Comparing the two sets of relations, we find that,
since $K$ and $\tilde K$ are not torsion classes,
the only possibility is
$\tilde K=(2r-1)K$ and $\tilde K=(2\tilde r-1)K$ with
 $r=\tilde r=1$. Therefore
 $K=\tilde K$.

\noindent{\it Case 2. $M$   rational or irrational ruled}.

\noindent In this case it is proved in [LL2]. For the
 convenience of readers,
we briefly present the argument here. Let us recall that a classical
theorem of Wu states that a class $c$ is the first Chern class of
an almost complex structure on $M$ only if $c^2=2\chi(M)+3\sigma(M)$
and its mod $2$ reduction is $w_2(M)$. This fact alone determines
the choice of the symplectic canonical classes up to sign
when $M$ is  $CP^2$ or an
$S^2-$bundle over $S^2$ via a simple calculation.

When $M$ is irrational ruled, the main observation is  that, because $M$ has
a metric of positive scalar curvature,
by the above mentioned result in [W] and Theorem 2.2, one concludes that the wall crossing number of $K^{-1}$
must be one.
By Lemma 3.3,
$$K^{-1}\cdot \gamma=2 SW^0_-({\cal L}_{K^{-1}})-2 SW^0_+({\cal L}_{K^{-1}})=2 \eqno (3.2).$$ Invoking Wu's theorem,
we again find the choice of $K$ is unique up to sign.

Now suppose $M$ is non-minimal. Let $N\# l\overline{CP}^2$
be a minimal reduction of $M$. Suppose $V$ is a symplectic canonical class
of $N$, then we have just proved that the only other
symplectic canonical class on $N$ is $-V$.
Given any symplectic canonical class $K$ of $M$, by Proposition 3.5,
there is a diffeomorphism
carrying $K$ to $V\pm F_1+...\pm F_l$ or $-V \pm F_1+...\pm F_l$,
Using reflections along $F_i$,
we see all the symplectic canonical classes
of $M$ can be carried to either $V+F_1+...+F_l$ or $-(V+F_1+...+F_l)$.
Therefore they are all equivalent. The proof of Theorem 1 is finished.

Let ${\cal K}$ be the set of
the symplectic canonical classes.
We can often  give a concrete description of
${\cal K}$.

When $M$ is rational  and $b^-\leq 9$
($b^-$ is the dimension of a maximal subspace of $H^2(M;{\bf R})$),
${\cal K}$ is the set of characteristic classes with square
$2\chi(M)+3\sigma(M)$,
since it was shown in [LL2]
that the group of  diffeomorphism acts
transitively on this set.

When $M$ is not rational nor ruled, let
$V$ be a symplectic canonical class of $N$, it is easy to see from
Theorem 1 and Proposition 3.6 that
$${\cal K}=\{\pm V\pm F_1\cdots \pm F_l\}.$$

\noindent {\bf Proposition 3.7}.  Suppose $M$ is irrational ruled
and is given a standard minimal reduction
$(S^2\times \Sigma_h)\# l\overline {CP}^2$
%Let  $F_i$ be
%the generators of  $H^2$ of the $\overline{CP}^2$.
and a standard basis.
 Let $\bf D$ be the set of  $(l+1)-$tuples
${\bf d}=\{\epsilon,c_1,...,c_l\}$ with $\epsilon=\pm2$ and $c_i$ odd.
For each $\bf d\in \bf D$ define
$$K_{\bf d}=\epsilon U +{(8-8g-l+\sum_i c_i^2)\over2\epsilon}T+   \sum c_i F_i.$$
Then ${\cal K}=\{K_{\bf d}|{\bf d}\in {\bf D}\}.$

\noindent{\it Proof}.
In this case $\gamma=T$. By the equation (3.2),
  if $K=aU+bT+\sum_i c_iF_i$ is in ${\cal K}$,
then $a=\pm2$. And since $K$ is characteristic with square $8-8g-l$,
$K$ must be of the form $K_{\bf d}$ for some ${\bf d}\in {\bf D}$.

To show that $K_{\bf d}\in{\cal K}$, it suffices to show
that $K_{0}$ can be carried to $K_{\bf d}$ by a
composition of reflections along classes in ${\cal E}$.
Since $T$ and $F_i$ are   represented by disjoint spheres,
and $T$ is of  square zero, it is easy to see that
$-kT-F_i$ is in ${\cal E}$ for any $k$ and $i$.
Let $e=aU+bT+\sum_i c_iF_i$.
Under the reflection $r^{i}_{-1}$ along $F_i$,
$$c_i\longrightarrow -c_i,  a\longrightarrow a, b\longrightarrow b, c_j\longrightarrow c_j \hbox{ if }j\ne i.$$
Under the reflection $f^k_i$ along $-kT-F_i$,
$$c_i\longrightarrow 2ka-c_i, b\longrightarrow b+2(-ka+c_i)k, a\longrightarrow
a,
c_j\longrightarrow c_j \hbox{ if }j\ne i.$$

Let $\bf d\in \bf D$ be a sequence with $\epsilon=2$.
 Write
$c_i$ as $4k_i-\tau_i$ with $\tau_i=\pm1$. Denote the identity on
$H^2(M;{\bf Z})$ also by $r^i_1$  for each $i$.
 Then
$$K_{\bf d}=(r_{\tau_1}^1\circ f_1^{k_1})\circ...\circ (r_{\tau_l}^l\circ
f_l^{k_l}) (K_0) \eqno (3.3)$$
and thus is in ${\cal K}$.
For the case $\epsilon=-2$, we just have to
observe that $K_{-\bf d}=-K_{\bf d}$.
In fact we have just proved that  the  set ${\cal K}$ is given by
$\{K_{\bf d}, \bf d\in \bf D\}$.

\noindent{\bf Corollary 1}. Let $M$ be a smooth, closed oriented
$4-$manifold. The number of equivalent classes of the symplectic canonical
classes is finite.

\noindent{\it Proof of Corollary 1}. If $M$ is a smooth, closed
oriented $4-$manifold with $b^+>1$,
the number of Spin$^c$ structures ${\cal L}$ such that $SW({\cal L})$
is non-trivial   is finite
([W]). This corollary simply
follows from this fact, Theorem 2.2 and Theorem 1.

\noindent{\bf Corollary 3.8}. Let $M$ be a minimal  $4-$manifold with
$b^+=1$ and non-empty symplectic cone. Suppose $K$ is a symplectic canonical class. If $\phi$ is a diffeomorphism, then $\phi^* K=\pm K$.

\noindent{\it Proof}. Let $\phi$ be a
diffeomorphism. Since
$\phi^*K$ is the symplectic canonical class of $\phi^*\omega$,
This corollary  is immediate from the proof of Theorem 1.

It follows from Corollary 3.8, up to sign, the action of $\phi$ on $H^2$ is
 determined by its
restriction to $L$, the orthogonal complement of $K$ in $H^2(M;{\bf Z})$.
If $K^2>0$, then $L$ is negative definite, and therefore
it has only finitely many
automorphisms. Since $K^2=2\chi(M)+s\sigma(M)>0$, we have the following generalization of
Corollary 4.8 in [FM2],

\noindent{\bf Corollary 3.9}. Let $M$ be a minimal closed, oriented
 $4-$manifold with
 $b^+=1$ and non-empty symplectic cone.
 Let $D(M)$ be the image of
diffeomorphisms of $M$ in the automorphisms of $H^2(M;{\bf Z})$. Then $D(M)$
is finite if $2\chi(M)+3\sigma(M)>0$.

\bigskip
\noindent{\bf \S4. ${\cal C}$ and ${\cal C}_K$}
\medskip

In this section, we are going to determine the
$K-$symplectic cone and the symplectic cone for
a closed oriented 4-manifold with $b^+=1$ and non-empty symplectic cone.

We start with some general properties of the $K-$symplectic
cone.

\noindent{\bf Proposition 4.1}. Let $M$ be a closed, oriented $4-$manifold
and $K$ be a symplectic canonical class.

\noindent 1. If $K'$ is another
symplectic canonical class, then
${\cal C}_{ K}\cap {\cal C}_{ K'}$ is empty.

\noindent 2. Let $f$ be a
diffeomorphism, then $f^* {\cal C}_{ K}={\cal C}_{ f^*K}$.

\noindent 3. Suppose $M$ has $b^+>1$ and $b_1,..., b_l$ are the
classes such that $SW({\cal L}_{K^{-1}}\otimes b_i)$ is non-trivial.
Then
$${\cal C}_{K}\subset \{e\in {\cal P}
|e\cdot K\geq 0 \hbox { and }
e\cdot K>|e\cdot b_i| \hbox {for all $b_i\ne \pm K$ }\}.$$

\noindent{\it Proof}. We first claim that $K-K'$
can not be a non-zero  torsion class. In the case $b^+>1$, this follows from
Theorem 2.4.1; in the case $b^+=1$, this is due to Theorem 1.
Thus if $e\in {\cal C}_K$ and $e'\in {\cal C}_{K'}$, we have by Theorem 2.4
$$K\cdot e>K'\cdot e  \hbox{     and     } K\cdot e' < K'\cdot e'.$$
Consequently $e\ne e'$ and the first part is proved.

 If $\omega$ is a symplectic form with
$K$ as its symplectic canonical class,  then
$f^*\omega$ is symplectic with $f^*K$ as its symplectic canonical class. Therefore
the second part  holds.

The last part follows directly from Theorem 2.4.
The proof of the proposition is complete.

We would like to remark that, for all the known manifolds with
trivial $K$, which is either the $K3$ surface,
or a $T^2-$bundle over $T^2$,  the third part is in fact an equality. Indeed if
 $K$ is trivial, then
the only Seiberg-Witten basic class is $0$ and
the righthand side is just ${\cal FP}$.  When $M$ is
a $T^2-$bundle over $T^2$, $K$ is trivial and it has been shown
explicitly in [G] that all classes in ${\cal P}$ can indeed be represented by
symplectic forms. For $K3$ surface, this is also the case.

\medskip
\noindent{\bf \S4.1. When $M$ is minimal}
\medskip

In this subsection, $M$ is a minimal closed, oriented $4-$manifold with
$b^+=1$.
We will first describe the set $A_K$. The knowledge of $A_K$
is then used to provide a complete description of ${\cal C}_K$.

\noindent{\bf Proposition 4.2}. Let $M$ be a minimal symplectic $4-$manifold
with
$b^+=1$. Let $e\in H^2(M)$ be a class in the forward cone. If $e-K$ is in
 the closure of the forward cone and is not equal to zero, then
$e$ is represented by connected symplectic surfaces.
In particular  for
$N$ big,  $Ne$ is represented by connected symplectic
surfaces.

\noindent{\it Proof}. The assumption that $e$  being in the
forward cone and $e-K$ in the closure of the forward cone implies
that $(e-K)\cdot e> 0$. Since $d(e)=(e-K)\cdot e$ is even, we have
$d(e)\geq 2$. It also implies that $2e-K=e+(e-K)$ is in the
forward cone, thus $(K^{-1}+2e)\cdot \omega>0$. By Lemmas
3.2$-$3.4, $SW^0_{\omega,-}({\cal L}_{K^{-1}}\otimes e)$ or
$SW^2_{\omega,-}({\cal L}_{K^{-1}}\otimes e)$ is nontrivial if $M$
is ruled or $d(e)\geq b_1$. But if $M$ is non-ruled, then
$b_1(M)\leq 2$ by [Liu]. Therefore, under the assumption,
$SW_{\omega,-}(K^{-1}\otimes e)$ is non-trivial. By Theorems
2.9(1) and 2.7, $e$ is represented by an embedded symplectic
surface.

Finally we prove this surface is connected. Since $M$ is assumed
to be minimal, every component has non-negative square.
Since $e$ is in the forward cone, $e^2$ is positive.
Therefore at least one component has positive square.
If the surface has more than one component,
then it violates the Light cone lemma 3.1.
 The proof is complete.

\noindent{\bf Theorem 2}. Let $M$ be a minimal closed, oriented $4-$manifold
with
$b^+=1$ and $K$ be a  symplectic canonical class.  Then
$${\cal C}_{M,K}={\cal FP}(K).$$
Consequently, any real cohomology class of positive square
is represented by an orientation-compatible symplectic form.

\noindent{\it Proof}. Fix a symplectic form $\omega$ whose symplectic
canonical class is $K$.
Since being symplectic is an open condition, we can assume
that $[\omega]$ is an integral class.

We first show that any integral $e$ in the forward cone
 is in the $K-$symplectic cone.
The first step is to show that, for a large integer $l$, $le-[\omega]$
is represented by a symplectic surface.
Indeed if $l$ is large, since $\omega\cdot e>0$,
 $le-[\omega]$ is in the forward cone and
$(le-[\omega])-K$ is in the closure of the forward cone. Thus by
Proposition 4.2, $le-[\omega]$ can be represented by a symplectic
surface. Given a symplectic surface $C$ with non-negative
self-intersection, the inflation process of Lalonde and McDuff
constructs a judicious Thom form $\rho$, representing the
Poinc\'are dual to the class of $C$ and supported in an
arbitrarily small neighborhood of $C$, such that
$\omega+\kappa\rho$ remains symplectic for all positive number
$\kappa$. Thus $le=[\omega] + (le-[\omega])$ is represented by a
symplectic form.

Since any positive real multiple of a symplectic form
is a symplectic form with the same canonical class,
 we have shown that any real multiple of an integral class
is in the $K-$symplectic cone.
 To show this is true  for
a general class $\alpha$ in the forward cone, we use  a trick in [Bi1]:
$\alpha$ can be written as $\alpha=\sum_{i=1}^p
 \alpha_i$, where the rays of $\alpha_i$ are arbitrarily close to
that of $\alpha$ and
each $\alpha_i=s_i\beta_i$
for some positive real number $s_i$
and  an integral class $\beta_i$.
Fix such a decomposition such that each $beta_i$ is in the forward cone.
Our strategy is to show inductively that for any $q$, $\sum_{i=1}^q \alpha_i$
is in the $K-$symplectic cone.

Since we know $\alpha_1$ is in the $K-$symplectic cone, we can
choose a symplectic form $\omega_1$ with $\alpha_1=[\omega_1]$.
For a large integer $l$, since $\beta_2\cdot \omega_1>0$,
by Proposition 4.2,
$l\beta_2$ is represented by
an $\omega_1-$symplectic surface.
By the inflation process, $\alpha_1+\kappa(l\beta_2)$ is
represented by a symplectic form for any real number $\kappa$.
If we choose $\kappa=s_i/l$, then we find that
$\alpha_1+\alpha_2$ is
in the $K-$symplectic cone.
Now choose a symplectic form $\omega_2$ with
$[\omega_2]=\alpha_1+\alpha_2$. It can be shown in the same way
that $l\beta_3$ is represented by a $\omega_2-$symplectic
surface and  $(\alpha_1+\alpha_2)+\alpha_3$ is
in the $K-$symplectic cone. Repeating this process,
we find that $\alpha$ is in the $K-$symplectic cone.
Thus we have shown that ${\cal C}_K={\cal FP}.$

The last statement is clear now because any nonzero real
multiple of
an orientation-compatible symplectic form is still
such a form. The proof is complete.

Let us remark that Donaldson's result mentioned in \S1
 fits nicely with our results.
Indeed, if $e\in {\cal FP}$,  Theorem 2 tell us that $e$ can be represented
by a symplectic form $\omega$. Donaldson's construction then produces
$\omega-$symplectic surfaces representing large multiples  of $e$.

\medskip
\noindent{\bf \S4.2. When $M$ is not minimal}
\medskip

In this subsection $M$ is a non-minimal 4-manifold
with $b^+=1$ and non-empty symplectic cone.

The following result is the analogue of Proposition 4.2.

\noindent{\bf Proposition 4.3}.
 Let $M$ be a symplectic $4-$manifold with
$b^+=1$ and symplectic canonical class $K$.
% Let ${\cal E}_{\omega}$ be the set of
%the exceptional classes represented by symplectic $-1$ spheres.
Let $e\in
H^2(M)$ be a class in the forward cone.
Assume that $e-K$ is in the closure of forward cone and $e-K\ne 0$. Further assume that
$e\cdot E\geq -1$ for all
$E\in {\cal E}_{\omega}$.
Then  $e$ can be represented by a symplectic surface.
Furthermore, if $e\cdot E\geq 0$ for all $E\in {\cal E}_{\omega}$, then the symplectic surface is
connected.

\noindent{\it Proof}.
Let $N$ be a (symplectic) minimal reduction of $M$, i.e. $N$ is minimal and $M$ is obtained
from $N$ by blowing up some number of points. Let $F_1, ..., E_k$
be the exceptional classes for the blow down map $M\longrightarrow N$
and $V$ be the symplectic canonical class of $N$.
Thus $F_i\in {\cal E}_{\omega}$ and $K=V+\sum F_i$.
Write $e$ as $e'-\sum n_i F_i$. Then
$$\eqalign {&e'\cdot e'=e\cdot e+\sum n_i^2 \cr
 &e'\cdot V=e\cdot K-\sum n_i\cr
&e'-V=(e-K)+\sum (n_i+1) F_i.\cr}$$
So $e'\cdot e'>0$ and $e'-K_N$ is in the closure of the forward cone.
If we use  $Gr^i(M,e)$ to denote the part of
$Gr_{\omega}(M,e)$ in $\Lambda^iH^1(M;{\bf Z})$, then, by Proposition 4.2,
  $Gr^0_{\omega}(N, e')$ or $Gr^2(N, e')$ is nonzero. Since $d(e)\geq 2$,
by Theorem 2.8, $Gr^0_{\omega}(M, e)$ or $Gr^2_{\omega}(M,e)$ is nonzero.
Now the proposition follows from Theorems 2.9(1) and 2.7.

A  slightly weaker version of Propositions 4.2-4.3  appeared in [MS] and [Mc3]
(see Lemma 2.2 in [Mc3] and the proof of Proposition 4.11 in this paper).
I. Smith informed us that via a construction using Lefschetz fibrations,
he and Donaldson  obtained similar  results about the existence of
symplectic surfaces.

For any two subsets  $U$ and $F$ of $H^2(M)$, we define
$$U_{F}=\{e\in U|e\cdot f>0 \hbox { for any $f\in F$ }\}.$$

\noindent{\bf Theorem 3}.    Let $M$ be a smooth, closed oriented
$4-$manifold with $b^+=1$ and $K$ be a symplectic canonical  class.
Use any class in $\Omega_K$ to define the forward cone.
Then
$${\cal C}_{K}={\cal FP}_{{\cal E}_K}$$

\noindent{\it Proof}.
By Lemma 3.5(1) ${\cal C}_K$ is contained in ${\cal FP}_{{\cal E}_K}$.
 To prove the inclusion in the other direction, we
fix any integral symplectic form $\omega$ on $M$ with $K$ as its symplectic canonical
class. Let $e$ be an integral class in ${\cal FP}_{{\cal E}_K}$. Since for
a large integer $l$,
$le-[\omega]$ and $(le-[\omega]-K$ are both  in the forward cone,
 and $le\cdot E>0$ for any $E\in {\cal E}_{\omega}$ by Lemma 3.5,
 $le-[\omega]$ is represented by a connected
$\omega-$symplectic surface by Proposition 4.3.
By the inflation process, we find that $e$ is in the $K-$symplectic
cone.

To pass to a general class $\alpha\in {\cal FP}_{{\cal E}_K}$, we notice that Biran's trick still
works here.  Since ${\cal FP}_{{\cal E}_K}$ is open,
$\alpha$ can be written as a sum $\sum_i \alpha_i$, where
$\alpha_i=s_i\beta_i$ for some positive real number and
$\beta_i$ is integral and lies in ${\cal FP}_{{\cal E}_K}$.
The rest of the proof is exactly the same as that in Theorem 2.
The proof is complete.

We would like to point out that an amusing fact from
Theorems 2 and 3 is that, in the case $b^+=1$,
${\cal C}_K$ is a convex set, which  is not obvious at
all from its definition.

We would also like to speculate on the relation between
the K\"ahler cone and the $K-$symplectic cone for K\"ahler surfaces
with $b^+=1$ (or equivalently zero geometric genus). q

In the case of
a rational surface, it is not known whether the K\"ahler cone can be as big as the symplectic cone.
By Nakai-Moishezon's criterion for ampleness, for a rational surface,
the K\"ahler cone is simply the set of real cohomology classes
with positive square,
which are positive on any irreducible curves with negative squares.
 In light of Theorem 5,
this is true if one can show that
for any integer $s$, there are $s$ generic points on the projective plane such that,
on the  rational surface obtained by  blowing  up these  points,
the only irreducible curves with negative squares are smooth
rational curves with square $-1$. It was observed in [MP] that this kind of
question is related to Nagata conjecture on the existence
of equisingular plane curves.
The connection was also pointed out to the authors by Koll\'ar.
In fact some very interesting results addressing the connection has been
obtained by Biran (see [Bi3]).
In [FM],   good generic rational surfaces are studied. For those
surfaces, the only irreducible
curves with negative squares
are the exceptional curves and those representing $K^{-1}$. Thus
by Nakai-Moishezon's criterion for ampleness, the K\"ahler cone is given by
$$\{e\in {\cal FP}|e\cdot C>0 \hbox{ for any holomorphic $-1$
 curve $C$ and $e\cdot K^{-1}>0$ }\}.$$

For general symplectic 4-manifold with $b^+=1$,
for an almost complex structure $J$ with $K$ as the canonical class,
we can introduce the $J-$symplectic cone
$${\cal C}_{ J}=\{e\in H^2(M;{\bf R})|e \hbox{ has a
symplectic representative
compatible with } J\}.$$
When $J$ is integrable, ${\cal C}_{ J}$ is just the K\"ahler cone.
It would be nice  to know when there
exists $J$ such that ${\cal C}_{ J}$=${\cal C}_{ K}$.

Now we will give prove Theorem 4,  which characterizes ${\cal C}$
in terms of the set ${\cal E}$. For this purpose we summarize some
facts about ${\cal E}$ in the next lemma.

\noindent{\bf Lemma 4.4}. Let $M$ be a closed oriented $4-$manifold
with ${\cal C}$ non-empty.
Suppose $M$ has a (smooth) minimal reduction $N\# l\overline{CP}^2$.
Let $F_1,...,F_l$ be the generators of $H^2$ of
 $\overline{CP}^2$.

\noindent 1. ${\cal E}=\cup_{K\in{\cal K}} {\cal E}_K$,
where ${\cal K}$ is the set of symplectic canonical classes.

\noindent 2.
If $e\in {\cal C}$ and $E\in {\cal E}$,
then $e\cdot E\ne 0$.

\noindent 3. When $N$ is not $CP^2$ nor an $S^2-$bundle,
${\cal E}=\{\pm F_1, ..., \pm F_l\}.$

\noindent{\it Proof}. Part 1 and 3 were  proved in [L1].
  We only have to prove part two.
When $b^+>1$, it is due to Taubes [T1].
The case $b^+=1$ is implicitly in [L1] and we will make it clear here.
Suppose $\omega$ is a symplectic form in the
class of $e$ and $K$ is its symplectic canonical class.
Since $K$ is characteristic and $E^2=-1$, $K\cdot E$ is odd. In particular
$K\cdot E=l \ne 0$.
Then by Lemma 2.2 and Theorem T1 in [L1], we find that $lE$
or $-lE$ is represented by a
symplectic surface (not necessarily embedded). This implies that $e\cdot E\ne 0$.
The proof is finished.

When $M$ is rational or irrational ruled with $b^-\geq 2$,
${\cal E}$ is in fact an infinite set. It is not hard to
write down it explicitly  when $M$ is irrational ruled.
When $M$ is rational,
it is determined   in [LiL2]
(it was also shown in [LL1] that  every class with square $-1$ is in ${\cal E}$
if $b^-\leq 9$).

We start the proof of  Theorem 4 in the easier case when
$M$ is not rational nor irrational ruled.

\noindent{\bf Proposition 4.5}.
 Let $M$ be a smooth, closed oriented $4-$manifold  with $b^+=1$ and
${\cal C}_M$ non-empty.
If $M$ is not rational nor irrational ruled, then
$${\cal C}_M=\{e\in {\cal P}|e\cdot E\ne 0 \hbox{ for any $E\in {\cal E}$}\}.$$

\noindent{\it Proof}.  Suppose $M$ is given a minimal reduction
$N\# l\overline{C P}^2$.
Let $F_1, ..., F_l$ be the generators of $H^2$ of the $\overline {CP}^2$.
By Lemma 4.4,  ${\cal E}=\{\pm F_1, ..., \pm F_l\}$.
If $e\in {\cal C}_M$, $e\in {\cal C}_{K}$ for some $K$.
By Lemma 3.5 and Theorem 3,
and $e\cdot F_i\ne 0$ for any $i$. Thus we have shown
$${\cal C}_M\subset \{e\in {\cal P}|e\cdot E\ne 0 \hbox{ for any $E\in {\cal E}$}\} .$$

Now let us prove the inclusion in the reverse direction.
Suppose  $e\in {\cal  P}$ and
$e\cdot F_i\ne 0$.
If $V$ is one of the two symplectic canonical classes of $N$,
as mention in \S3, any symplectic canonical
class of $M$ is of the form $K=\pm V\pm F_1\pm...\pm F_l$.
By possibly changing $e$ to $-e$, we can assume that $e$ is in the forward cone determined by any symplectic canonical class of the form
$V+\pm F_1\pm...\pm F_l$.
Let $\epsilon_i=(e\cdot F_i)/|e\cdot F_i|$ and  $K=V+\sum \epsilon_i F_i$.
Then ${\cal E}_K=\{\epsilon_1F_1,...,\epsilon_l F_l\}$,
and  $e\cdot \epsilon_i F_i>0$. Therefore  $e\in {\cal C}_{ K}$
by Theorem 3. Thus we have shown that $$\{e\in {\cal P}|e\cdot E\ne 0 \hbox{ for any $E\in {\cal E}$}\}\subset {\cal C}_M.$$
The proof of Proposition
4.5
is complete.

% When $N$ is $CP^2\#\overline{CP}^2$,
%${\cal E}=\{\pm F_1\}$.
%The same argument works for $CP^2\#\overline{CP}^2$ and the proof of

To deal with the non-minimal rational and irrational ruled 4-manifolds,
we need to introduce the notion of a reduced class.

\noindent{\bf Definition 4.6}.
For a non-minimal rational manifold with a standard
decomposition $CP^2\# l\overline{CP}^2$
and a standard basis $\{H, F_1,\dots, F_l\}$,
a  class  $\xi=aH-\sum^n_{i=1}b_iF_i$ is called reduced
if $$\cases{b_1\geq b_2\geq \cdots \geq b_l\geq 0\cr
a\geq                b_1+b_2+b_3.\cr}$$
For a non-minimal irrational ruled manifold with a standard decomposition
$S^2\times \Sigma_g\# l\overline{CP}^2$ and a
standard basis $\{U, T,F_1, \dots, F_l\}$, a class $e=aU+bT-\sum c_iF_i$ is called
reduced if
$$\cases{c_1\geq c_2\geq \cdots\geq c_l\geq 0\cr
a\geq c_i
\hbox{ for any }i.\cr}$$

The following facts about reduced classes in [LiL2] are crucial.

\noindent{\bf Lemma 4.7}. Let $M$ be a non-minimal
rational or irrational ruled
$4-$manifold with a standard decomposition and a standard basis.

\noindent 1. There is a simple algorithm to
transform any class of positive square   to a reduced class
via a  diffeomorphism.

\noindent 2.  Suppose $e$ is a reduced class and
$E\in {\cal E}_{K_0}$.  Then
 $e\cdot E\geq 0$, and $e\cdot E>0$ if $E$ is not one
of the $F_i$.

\noindent{\bf Proposition 4.8}. Let $M$ be a non-minimal rational
or irrational ruled 4-manifold.
Then $$ \{e\in {\cal P}|0< |e\cdot E|
 \hbox{ for any $E\in {\cal E}$}\}={\cal C}_M.$$

\noindent{\it Proof}.  Let $M$ be given a standard decomposition
and standard basis, so we can define reduced classes.
Let $e$ be a class in ${\cal P}$
such that $e\cdot E\ne 0$ for any $E\in {\cal E}$.
By Lemma 4.7(1) $e$ can be transformed to a reduced class
by some  diffeomorphism, say $\phi$.
On one hand, since the set ${\cal E}$ is preserved by the group of
diffeomorphisms,
$\phi^*e\cdot E\ne 0$ for any $E\in {\cal E}$.
On the other hand, since $\phi^*e$ is a reduced class, by Lemma 4.7(2),
$\phi^*e\cdot E\geq 0$ for any $E\in {\cal E}_{K_{0}}$.
Therefore $\phi^*e\cdot E> 0$ for any $E\in {\cal E}_{K_{0}}$.
By Theorem 3, we  see that $\phi^*e\in {\cal C}_{K_{0}}.$ Therefore
by Proposition 4.1, $e\in   \phi^*{\cal C}_{K_0}={\cal C}_{\phi^*K_0}\subset {\cal C}.$   Thus we have proved $$ \{e\in {\cal P}|0< |e\cdot E|
 \hbox{ for any $E\in {\cal E}$}\}\subset {\cal C}_M.$$

The inclusion in the reverse direction follows from Lemma 4.4 (2).
The proof of Proposition 4.8 is complete.

Propositions 4.5 and 4.8 immediately give us

\noindent{\bf Theorem 4}.
 Let $M$ be a smooth, closed oriented $4-$manifold  with $b^+=1$ and
${\cal C}_M$ non-empty.
Then
$${\cal C}_M=\{e\in {\cal P}|e\cdot E\ne 0 \hbox{ for any $E\in {\cal E}$}\}.$$

When $M$ is given a minimal reduction, we can
give a more explicit presentation of
${\cal C}$.

\noindent{\bf Proposition 4.9}. Suppose $M$ is as in Theorem 4
and is given a minimal reduction $N\# l\overline {CP}^2$.
Let  $F_i$ be
the generators of  $H^2$ of the
$\overline{CP}^2$.
When $M$ is rational or irrational
ruled, further assume the minimal reduction is  standard
and a standard basis is given.

\noindent 1. When $M$  is not rational nor irrational ruled,
A class $e$ with positive square is represented by a symplectic form
if and only if $e\cdot F_i\ne 0$ for any $i$.

\noindent 2. When $M$ is irrational ruled,
 a class $e=aU+bT+\sum_i c_iF_i$
with positive square is represented by a symplectic form if and
only if $C_i/A$ is not an integer for any $i$.

\noindent 3. When $M$ is rational, a reduced  class $e=aH-\sum^n_{i=1}b_iF_i$
with positive square
is in ${\cal C}_{K_0}$ if and only
if $b_i>0$ for each $i$. For a general class $e$,
transform it to a reduced class $e'$. Then $e$ is represented
by a symplectic form if and only if $e'$ is thus represented.

\noindent{\it Proof}. The first part is contained in the
proof of Proposition 4.5.
If $M$ is irrational ruled, the conclusion follows
immediately from

\noindent{\bf Lemma 4.10}.
If
$Y=\{s_1T\pm F_1,...,s_lT \pm F_l, s_i\in{\bf Z}\}$,
then ${\cal E}=Y$.

\noindent{\it Proof}. We have shown in course of
the proof of Proposition 3.7 that
$Y\subset {\cal E}.$
To prove the inclusion in the reverse direction, by Lemma
4.4(1), it suffices to show that for every symplectic canonical
class $K$,
${\cal E}_K\subset Y.$

Consider $K_0$ first.
If the  symplectic form $\omega$ is obtained by
blowing up a product symplectic form on $S^2\times \Sigma_h$
which
 pairs positively with $U$ and $T$,
then it has $K_0$ as the canonical class.
For such a form, it is not hard to show (see [Bi1])
that
$${\cal E}_{K_0}=\{F_1,..., F_l, T-F_1,..., T-F_l\}.$$
So ${\cal E}_{K_0}\subset Y.$
By Proposition 3.7, each symplectic canonical class is of the
form $K_{\bf d}$.
By equation (3.3) and Lemma 3.5,
$$\eqalign{{\cal E}_{K_{\bf d}}=&(r_{\tau_1}^1\circ f_1^{k_1})
\circ...\circ (r_{\tau_l}^l\circ
f_l^{k_l}) {\cal E}_{K_0}\cr
=&\{-2k_1T +\tau_1 F_1, (1+2k_1)T-\tau_1 F_1,
...,- 2k_lT+\tau_l F_l, (1+2k_l)T-\tau_l F_l                               \}. \cr}$$
Notice that $\tau_i=\pm 1$ for each $i$,
 ${\cal E}_{K_{\bf d}}\subset Y$. Lemma 4.10 is proved.

Suppose $M$ is rational.
Let  $e=aH-\sum^l_{i=1}b_iF_i$ be a reduced class
with positive square. If it is in ${\cal C}_{K_0}$,
then since $F_i\in {\cal E}_{K_0}$, it is necessary that
  $b_i>0$ for each $i$. Conversely, if  $b_i>0$ for each $i$,
by Lemma 4.7(2) and Theorem 3, it is in ${\cal C}_{K_0}$.
The first statement is thus proved. The last statement is obvious,
so the proof of Proposition 4.9 is finished.

As we previously remarked, when $M$ is rational,  though
${\cal E}$ can be determined, it is hard to
write down it explicitly. So we do not have a as nice presentation
of ${\cal C}$ as those in the other cases. However,
in light of Lemma 4.7(1), it is still a very effective one.

Having determined the image of the map $CC:\Omega_M\longrightarrow
H^2(M;{\bf R})$, we are  also able to
say something about its inverse image by generalizing a result
of McDuff.

\noindent{\bf Proposition  4.11}.    Let $M$ be  closed oriented
$4-$manifold with $b^+=1$. Let $\omega_1$ and $\omega_2$ be
two  cohomologous symplectic forms.
If they can be joined by a path of symplectic forms then
they can be joined by a path of cohomologous symplectic forms.

\noindent{\it Proof}. This result was  proved as Theorem 1.2
in [Mc3] under the assumption that $M$ is not  of `Seiberg-Witten simple type'.
In [Mc3], a symplectic 4-manifold is said to be
of `SW simple type' if its only nonzero $Gr^0$ invariant
occur in classes with zero Gromov-Taubes dimension.
For a symplectic 4-manifold $M$ not of `SW simple type',
what is used in the proof  in [Mc3] is the following fact:
assuming $[\omega]$ is rational, then
there is a basis of $H^2(M;{\bf Q})$ formed by rational
classes $n[\omega], e_1,..., e_k$ with
$e_j^2<0$ for all $j$ such that $l(n[\omega]\pm e_j)$ is
represented by a connected $\omega-$symplectic surface
for all $j$ and large $l$.

By  Propositions 4.1-4.3,
all symplectic $4-$manifolds with $b^+=1$
satisfy this property. So the proof of this Proposition is identical
to that   of Theorem 1.2 in [Mc3].

In terms of the map $CC$, the result above can be interpreted as saying that,
when restricted to a connected component of $\Omega_M$,
the inverse image of a point of $CC$ is connected.
  This result
is  useful in \S6 for the duality conjecture.

\medskip
\noindent{\bf \S.5. Symplectic Casteluovo's criterion}

\medskip
In this and the following sections we will study further $A_K$, the
set of `$K-$stable' classes of symplectic surfaces.
Here we will determine, for a minimal 4-manifold,  which
integral multiples of $K$
is in $A_K$.  As an interesting corollary, we obtain
the sympletic Casteluovo's criterion for rationality.

\noindent{\bf Lemma 5.1}.  Let $M$ be a minimal  4-manifold
with $b^+=1$ and $K$ a non-torsion symplectic canonical
class. $A_K$ contains only non-torsion classes
with non-negative square.

\noindent{\it Proof}. We observe that, by definition, if $a\in A_K$ and
$e\in {\cal C}_K$, then $a\cdot e>0$. So obviously $a$ can not be a torsion
class.
Suppose $a\in A_K$  has negative square.
The orthogonal complement of $a$ in $H^2(M;{\bf R})$
still contains classes of positive square.
Let  $\beta$ be such a class.
 By Theorem 2, either $\beta\in {\cal C}_K$
or $-\beta\in {\cal C}_K$. However, the fact that
$a\cdot \beta=0$ contradicts with
the observation above. The proof is finished.

Notice  that,
for a fixed symplectic structure, there might be classes of negative square
 represented by symplectic surfaces. Such examples are easy to find.
For any positive integer $n$, the Hirzebruch surface
$F_2n$ is a minimal algebraic surface having holomorphic curves (hence symplectic with respect to
any K\'ahler form)
of square $-2n$ in  $F_{2n}$.

\noindent{\bf Proposition 5.2}. Let $M$ be a minimal  4-manifold
with $b^+=1$ and $K$ a non-torsion symplectic canonical
class. Then  $nK$ is
in $A_K$ in the following case:

\noindent 1. $n\leq -1$ and $M$ is $CP^2$, $S^2\times S^2$ or an $S^2-$bundle over $T^2$,

\noindent 2. $n=1$  and  $M$ has $b_1=2$ and is not an $S^2-$bundle over $T^2$,

\noindent 3.   $n\geq 2$ and $M$ is not rational nor irrational ruled.

%\noindent In fact in case 1 and 2,
%$Gr_{\omega}(nK)$ is non-trivial.

%\noindent 1. If $K$ is not in ${ S}_K$ and $b_1\geq 2$, then $M$ is ruled.
% If $K$ is not in ${ S}_K$ then $M$ is rational and ruled.

\noindent{\it Proof}.
By Lemma 5.1, for any pair $M,K$ with $K^2<0$,  no multiple of $K$ is in $A_K$.
Such manifolds
 are $S^2-$bundles over Riemann surfaces of genera at least 2.

Assume now that $M$ is not an $S^2-$bundle over a Riemann surface of genus at least 2.
 Then for any integer $n$, the Seiberg-Witten dimension of ${\cal L}_{K^{-1}}\otimes nK$ is
non-negative, because
 $$d({\cal L}_{K^{-1}}\otimes nK)=-K\cdot nK + nK \cdot nK= (n^2-n)K^2\geq 0.$$
Let $\omega$ be a symplectic form with  $K$ as its symplectic canonical
class.

The manifolds in case 1 satisfy
$K^2\geq 0$, $K\cdot \omega<0$.
So  $K^{-1}+2nK$ is in the forward cone for all $n\leq -1$.
Thus it follows from Lemma 3.4 that $Gr_{\omega}(nK)$ is non-trivial for $n\leq -1$
and case 1 is settled.

%Since other manifolds all satisfy $K\cdot \omega\geq 0$, this case is settled.

The manifolds in case 2 satisfy $K^2\geq 0$ and $K\cdot \omega>0$.
Let $\gamma$ be the generator
of $\Lambda^2H^1(M;{\bf Z})$ such that $\omega\cdot \gamma\geq 0$.
By Lemma 3.1, $K\cdot \gamma\geq 0$.
By Theorem 2.2 and Lemma 2.3 $SW^0_{\omega,+}({\cal L}_{K^{-1}}\otimes
K)=1$.
When $b_1=2$, by Lemma 3.3,
 $$SW^0_{\omega,-}(K^{-1}\otimes K)=-SW^0_{\omega,+}(K^{-1}\otimes K)
+{1\over 2}K\cdot \gamma\ne 0$$
unless  $K\cdot \gamma=-2$. But this is impossible, so by
Theorem 2.9 $Gr_{\omega}(K)$ is non-trivial.
Case 2 is settled.

By Lemma 3.2, $SW_{\omega,+}({\cal L}_{K^{-1}}\otimes nK)$ is zero if $n\geq 2$.
Since the manifolds in case 3 satisfy   $K\cdot \omega>0$, $(K^{-1}+2nK)\cdot
\omega>0$ for $n\geq 1$.
When $b_1=0$, or when $b_1=2$ and $K\cdot \gamma\ne 0$, by Lemma 3.3,
 $SW^0_{\omega,-}(K^{-1}\otimes K)-SW^0_{\omega,+}(K^{-1}\otimes K)$
is nonzero. Therefore
$Gr_{\omega}( nK)$ is nonzero for $n\geq 2$.

When $b_1=2$ and $K\cdot \gamma=0$, from part 2 we know $Gr_{\omega}(K)$
is non-trivial. By Lemma 3.1, $K$ must have square zero, so it is  represented by an embedded
symplectic torus. Thus, for each $n\geq 2$,  $nK$ can be represented by an embedded symplectic torus as well by Lemma 2.6.
 The proof of Proposition 5.2 is complete.

Now we can give the proof of the symplectic Castlenuova criterion.

\noindent{\bf Corollary 2}. Let $M$ be a closed, oriented $4-$manifold with
$b^+=1$ and a non-torsion symplectic canonical class  $K$.
If $b_1=0$ and $2K$ is not in $A_K$, then $M$ is rational.

\noindent{\it Proof}. First assume that $M$ is minimal. By Corollary 5.2, $M$ must be rational or irrationally ruled.
Since irrational ruled manifolds have $b_1>0$, $M$ is rational as claimed.

\medskip
\noindent{\bf \S6. $K-$surface cone and duality conjecture}
\medskip

Recall that the rational $K-$surface cone ${\cal S}^{\bf Q}_K$, introduced in \S1, is the cone
$\sum_{v\in A_{K}}{\bf Q}^{+}v$ in $H^2(M;{\bf Q})$.
In this section we will study this cone and
discuss the duality between it and the rational $K-$symplectic cone ${\cal C}^{\bf Q}_K={\cal C}_K\cap H^2(M;{\bf Q})$
inside $H^2(M;{\bf Q})$.

We start with a couple of simple algebraic lemmas.
 Let $W$ be a subset in $H^2(M;{\bf Q})$.
Define the dual of $W$ to be
$$W^{\wedge}=\{v\in H^2(M;{\bf Q})| v\cdot w>0 \hbox{ for any $w\in W$} \}. $$
Clearly $W^{\wedge}$ is a convex subset, and $W^{\wedge \wedge}=W$
if $W$ is convex.

\noindent{\bf Lemma 6.1}. Let $M$ be a symplectic $4-$manifold with
symplectic canonical class $K$.
Then ${\cal S}^{\bf Q}_K$ is convex, and
 $$ {\cal S}^{\bf Q}_{ K}\subset  {{\cal C}^{\bf Q}}^{\wedge}_ K \hbox{       and      }
{{\cal C}^{\bf Q}_{ K}}\subset {{\cal S}^{\bf Q}_K}^{\wedge}.$$

\noindent{\it Proof}. It directly follows from the definitions
of the two cones.

\noindent{\bf Lemma 6.2}.  Let $M$ be a closed oriented $4-$manifold
with $b^+=1$. Let ${\cal FP}^{\bf Q}$ be a component of ${\cal P}^{\bf Q}$ and $F$ be
a subset of $H^2(M;{\bf Q})$.
Then $${{\cal FP}^{\bf Q}_F }^{\wedge}=\overline {\cal FP}^{\bf Q}+\sum_{f\in F} {\bf Q}^+ f.$$
%And the same holds if ${\cal FP}^{\bf Q}$ is replaced by ${\cal P}$.

\noindent{\it Proof}. Let us first show that $\overline {\cal FP}^{\bf Q}=
{{\cal FP}^{\bf Q}}^{\wedge}$.  By the light cone lemma, ${\cal FP}^{\bf Q}\subset
{{\cal FP}^{\bf Q}}^{\wedge}$.  By the argument in Lemma 5.1 ${{\cal FP}^{\bf Q}}^{\wedge}\subset \overline
{\cal FP}^{\bf Q}$.
Since
$$\sum_{f\in F} {\bf Q}^+ f\subset {{\cal FP}^{\bf Q}_F}^{\wedge}
\hbox {     and     } \overline{\cal FP}^{\bf Q}={{\cal FP}^{\bf Q}}^{\wedge}\subset
{{\cal FP}^{\bf Q}_F}^{\wedge},$$
 we have, by the convexity of ${{\cal FP}^{\bf Q}_F}^{\wedge}$,
$$\sum_{f\in F} {\bf Q}^+ f + \overline{\cal FP}^{\bf Q}\subset
{{\cal FP}^{\bf Q}}^{\wedge}_F.$$
Let  $p\in \overline {\cal FP}^{\bf Q}+\sum _{f\in F} {\bf Q}^+ f$.
Then $p\cdot e>0$ for any $e\in \overline {\cal FP}^{\bf Q}$, thus
$p\in \overline {\cal FP}^{{\bf Q}^{\wedge}}={\cal FP}^{\bf Q}.$
And $p\cdot f>0$ for any $f\in F.$ Therefore
$$(\overline  {\cal FP}^{\bf Q}+\sum _{f\in F} {\bf Q}^+ f )^{\wedge}\subset {\cal FP}^{\bf Q}_F.$$
The proof is complete.

\noindent{\bf Theorem 5}.  Let $M$ be a smooth, closed oriented
$4-$manifold with $b^+=1$. Let $K$ be a symplectic canonical  class.
Use any class in $\Omega_K$ to define the forward cone.
Then $${\cal FP}^{\bf Q}+
\sum_{E_i\in {\cal E}_K} {\bf Q}^+ E_i \subset  {\cal S}^{\bf Q}_{ K}
 \subset \overline {\cal FP}^{\bf Q}  +
\sum_{E_i\in {\cal E}_K} {\bf Q}^+ E_i.$$

\noindent{\it Proof}.
%By Proposition 4.3, $ {\cal S}^{\bf Q}_K$
%contains the set  $$\{e\in {\cal FP}^{\bf Q}|e\cdot E_i\geq 0 \hbox { for all
%$E_i\in {\cal E}_K$ } \}.$$
We start with the first inclusion.
Let $\omega$ be a symplectic form with $K$ as its symplectic canonical
class.
By the same argument as in Propositions 4.2-4.3, we can show
that, if $e$ is an integral class in $ {\cal FP}^{\bf Q}$, then for a large
integer $l$, $SW_{\omega,-}({\cal L}_{K^{-1}}\otimes le)$ is non-trivial.
Therefore $le\in{\cal S}^{\bf Q}_{\omega}$ by Theorem 2.9(1).
We thus have shown that ${\cal FP}^{\bf Q}\subset {\cal S}^{\bf Q}_K$.
By Lemma 3.5, ${\cal E}_K={\cal E}_{\omega}$,  so
$$\sum_{E_i\in {\cal E}_K} {\bf Q}^+ E_i \subset  {\cal S}^{\bf Q}_{ K}.$$
By the convexity of ${\cal S}^{\bf Q}_K$, we have
$${\cal FP}^{\bf Q}+
\sum_{E_i\in {\cal E}_K} {\bf Q}^+ E_i \subset  {\cal S}^{\bf Q}_{ K}.$$

The second inclusion follows from Theorem 3 and Lemmas 6.1-6.2.
The proof is complete.

%Suppose $e\in {\cal FP}^{\bf Q}$. Define a subset $B$ of ${\cal E}_K$,
%$$B=\{E\in {\cal E}_K|E\cdot e<0.$$
%If suffices to show that there exists $e'\in

\noindent{\bf Duality Conjecture}. Let $M$ be a closed, oriented 4-manifold
with $b^+=1$.
Suppose $K$ is a symplectic canonical class, then
 the rational $K-$surface cone and the
rational $K-$symplectic cone are dual to each other.

We remark that this conjecture
 can be viewed as the analogue of the duality between the K\"ahler
cone and the cone of numerically effective curves on an algebraic
surface.

We will now show the Duality Conjecture
holds for several classes of minimal 4-manifolds.

\noindent {\bf Lemma 6.3}. Let $M$ be a closed oriented
4-manifold with $b^+=1$ and a symplectic canonical class $K$.
The conjecture holds for $M$ and $K$ if  any class $e \in \overline {\cal FP}^{\bf Q}(K)$ with square 0 is in $ {\cal S}^{\bf Q}_{K}$.

\noindent{\it Proof}. ${\cal C}^{\bf Q}_K$ is convex by Theorem 3, so ${{\cal C}^{\bf Q}_K}^{\wedge\wedge}={\cal C}^{\bf Q}_K$.
Therefore ${\cal C}^{\bf Q}_K={{\cal S}^{\bf Q}_K}^{\wedge}$ if ${\cal S}^{\bf Q}_K={{\cal C}^{\bf Q}_K}^{\wedge}$.
By Lemmas 6.1-2 and Theorems 3 and 5,
$${\cal FP}^{\bf Q}+
\sum_{E_i\in {\cal E}_K} {\bf Q}^+ E_i \subset {\cal S}^{\bf Q}_K\subset {{\cal C}^{\bf Q}_K}^{\wedge}={{\cal FP}^{\bf Q}_{{\cal E}_K}}^{\wedge}
=\overline  {\cal FP}^{\bf Q}  +
\sum_{E_i\in {\cal E}_K} {\bf Q}^+ E_i .$$
So
 ${\cal S}^{\bf Q}_K={{\cal C}^{\bf Q}_K}^{\wedge}$ if the complement of ${\cal FP}^{\bf Q}(K)$ in $\overline {\cal FP}^{\bf Q}(K)$ is in ${\cal S}^{\bf Q}_{K}$.
Since this complement consists of the classes $e \in \overline {\cal FP}^{\bf Q}(K)$ with square 0, the lemma is proved.

Two cohomologous symplectic forms are said to be isotopic if they can be
joined by a path of cohomologous symplectic forms. For a symplectic form
$\omega$, recall that $A_{\omega}$ is the set of integral classes which
can be represented by $\omega-$symplectic surfaces.

\noindent{\bf Lemma 6.4}. Let $M$ be a closed, oriented $4-$manifold
with $b^+=1$.
 If $\omega_1$ and $\omega_2$ are two symplectic forms which are deformation
equivalent and cohomologous, then $A_{\omega_1}=
A_{\omega_2}$.

% for isotopic symplectic forms $\omega_1$ and $\omega_2$, a class
%is represented by $\omega_1-$symplectic surfaces if and only if it
%is represented by $\omega_2-$symplectic surfaces, i.e.

\noindent{\it Proof}. By the basic result of Moser [M] it is easy to see that
$A_{\omega_1}=
A_{\omega_2} \hbox{ if $\omega_1$ and $\omega_2$ are isotopic. } $
Now the lemma follows from Proposition 4.11.

\noindent{\bf Proposition 6.5}.
The Duality conjecture holds for

\noindent a). 4-manifolds with  torsion symplectic canonical classes, $b^+=1$ and $b_1=0$.

\noindent b). $CP^2$ and $S^2\times S^2$.

\noindent c). $S^2-$bundles over $T^2$.

\noindent d).  $S^2\times \Sigma$ with $\Sigma$ a surface with genus at least two.

\noindent e). $T^2-$bundles over $T^2$ with $b^+=1$, $\Lambda^2H^2(M;{\bf Z})$ non-trivial,
and having a unique deformation
class of symplectic forms.

\noindent f). $S^1\times X$ with  $X$  a fibered $3-$manifold with
$b_1=1$ and the genus of the fiber at least 2, and having  a unique deformation
class of symplectic forms.

%\noindent f).

\noindent{\it Proof}.
For the manifolds in e) and f),  symplectic structures can be constructed by the construction in [Th].
So all the manifolds listed have symplectic structures.
Let $K$ be any one of the two symplectic canonical classes. Let $\omega$ be an arbitrary symplectic form
with $K$ as its symplectic canonical class.
By Lemma 6.3 what we need to  show is, for
any class $e$ in $\overline {\cal FP}^{\bf Q}(K)$ with $e^2=0$, some (rational) multiple of it is represented by an $\omega-$symplectic
surface.
We will show it is the case for manifolds in a)$-$e).

For manifolds in a), $K$ is a torsion class, so any class
$e$ in $\overline {\cal FP}^{\bf Q}(K)$ with square zero has $d(e)=0$ and
$(K^{-1}+2e)\cdot \omega >0$. If, in addition $b_1=0$, by Lemmas 3.1 and 3.2(1),
$SW_{\omega,-}({\cal L}_{K^{-1}}\otimes e)$ is non-trivial. By Theorem 2.9, $e\in A_{K}$.

For manifolds in b),
$K^{-1}$ is in the forward cone. Thus
$K^{-1}+2e$ is in the forward cone if $e$ is in $\overline {\cal FP}^{\bf Q}(K)$. By Lemma 3.4 and Theorem 2.9
$Gr_{\omega}(e)\ne 0$
for any class  $e\in \overline {\cal FP}^{\bf Q}$.

For manifolds in c), $K^{-1}$ has square zero and is in $\overline {\cal FP}^{\bf Q}(K)$. By the same argument as above,
if $e$ is in $\overline {\cal FP}^{\bf Q}(K)$ with square zero and $e\ne K^{-1}$, then it   is in $A_K$. The fact that
$K^{-1}$ is also in $A_K$ is proved in Proposition 5.3.

For a product irrational ruled 4-manifold $S^2\times \Sigma$,
 if $e$ has square $0$, then it is
Poincar\'e dual to  either a multiple of  $[S^2]$ or $[\Sigma]$.
On such a 4-manifold every symplectic form is isotopic to
a product from by [LM].
Since any product symplectic form remains a symplectic form on any $S^2$ and $\Sigma$, the conclusion
 follows from Lemma 6.4.

%For a non-product ruled surface,
%the classes with square zero are Poincar\'e dual to multiples of  $[S^2]$
%and $s=[s^+]+[s^-]$, where $s$ and $s'$ are any pair of sections
%with opposite squares. On such surfaces, every symplectic form is isotopic
%to a K\"ahler form. And for any integrable complex structure,
%there are a pair of holomorphic sections, which give
%rise to a symplectic surface representing $s$. The result follows from
%Corollary 4.12.

For the manifolds in e), $K$ is trivial and both $b_1$ and $b_2$
are equal to two. If $\gamma$ is the non-trivial generator of
$\Lambda^2H^2(M;{\bf Z})$ as in Lemma 3.3, it is the class of the
fibers. Since the intersection form on $H^2(M;{\bf Z}$ is
$(1)\oplus (-1)$, the classes with square zero are (rational)
multiples of $\gamma$ and another integral class which we denote
by $\eta$. We can assume that $\eta\cdot \omega>0$. Since $K=0$,
$d(\eta)$ is equal to zero. Notice that $\gamma\cdot \eta\ne 0$.
This, together with Lemmas 3.2-3, imply  that $SW_{\omega,-}({\cal
L}_{K^{-1}}\otimes \eta)$ is non-trivial. Therefore
$Gr_{\omega}(\eta)$ is non-trivial. All these manifolds are
geometric, and it is shown in [G] that every class with positive
square can be represented by a `geometric' symplectic form  such
that is symplectic on the fibers. In particular, the class of
$\omega$ is represented by such a `geometric' symplectic form
$\omega'$.
 Under the assumption that there is a unique
deformation class of symplectic forms, the conclusion follows from
Lemma 6.4.

For the manifolds in f),  the K\"unneth formula tells us that $b_1$ and $b_2$ are
again equal to two. Since $K^2=2\chi(M)+3\sigma(M)$, $K$ has square 0.
Let $g$ be the genus of the fibers of $\pi:X\longrightarrow S^1$.
Then $M=S^1\times X$ fibers over the 2-torus with fibers of genus $g$.
Since $g$ is at least 2, the class of the fibers is non-trivial in $H^2(M;{\bf Q})$.
In fact by the construction in [Th], the fibers are symplectic with respect to
some symplectic structures on $M$.
By the ajunction formula in [LL2], $K\cdot \gamma=(2g-2)$. Therefore $K$ is non-trivial in
$H^2(M;{\bf Q})$ as well and
the classes with square zero are
(rational) multiples of
$K$ and $\gamma$. By Proposition 5.2 $K$ is in $A_K$, so we only have to deal with
$\gamma$ by explicitly constructing a
symplectic form $\omega'$ in the class of $\omega$ such that
the fibers are $\omega'-$symplectic. Then, under the assumption that there is a unique deformation
class of symplectic forms, the conclusion that $\gamma$ is represented
by an $\omega-$symplectic surface again follows from
Lemma 6.4.
Choose a metric on $X$ such that $\pi:X\longrightarrow S^1$
is a harmonic map. Let $*_X$ be the star operator on $X$ with respect to this metric.
 Let $d\theta$ be the volume form of
the base circle, and $\mu$ be $*_X\pi^*d\theta$, . Let $ds$ be the
volume form of the product circle. Then for any pair of positive
numbers $\alpha$ and $\beta$, the form $\pm (\alpha ds\wedge
d\theta + \beta \mu)$ is symplectic and restricts to a symplectic
form on each fiber.  In particular the class of $\omega$ can be
thus represented. The proof of Proposition 6.5 is complete.

\noindent{\bf Remark 6.6}. It seems to us that it should not be too hard to show that
the duality conjecture holds for a manifold $M$
if it holds for one of its minimal models.
Given any symplectic form $\omega$ with $K$ as the symplectic
canonical class.
Blow down $M, \omega$ to a symplectic 4-manifold $N, \beta$.
Let $F_1,..., F_l$ be the
exceptional classes, and $V$ be the symplectic canonical
class of $\beta$.
If the conjecture holds for $N$, then
$\overline{\cal FP}^{\bf Q}_N(V)={\cal S}^{\bf Q}_V\subset {\cal S}^{\bf Q}_{\beta}$.
It is easy to see that
$$\overline {\cal FP}^{\bf Q}_M(K)\subset \overline{\cal FP}^{\bf Q}_N(V)+\sum_{i=1}^l {\bf Q}^+
F_i.$$
Thus it suffices to show that ${\cal S}^{\bf Q}_{\beta}=\overline {\cal FP}^{\bf Q}_N(V)\subset
{\cal S}^{\bf Q}_{K}.$
Due to Lemma 6.4,  this is quite possible.

Finally let us mention that we can also introduce the real
$K-$surface cone ${\cal S}_K$. It is easy to check that the
statements in Lemmas 6.1-3 and Theorem 5 remain valid for ${\cal
S}_K$ if the ${\bf Q}$ are removed or replaced by ${\bf R}$. The
real $K-$surface cone is similar to (often smaller than) the cone
of numerically effective curves. In fact, a more precise analogue
of the real $K-$surface cone is the deformed symplectic effective
cone introduced by Ruan in [R2]. In particular, if $\Omega_{K}$
has only one connected component, then they coincide.

\bigskip

\bigskip
\noindent{\bf Reference}
\medskip

\noindent [B] R. Brussee, The canonical class and the $C\sp \infty$
properties of KaŠhler surfaces. New York J. Math. 2 (1996), 103--146.

\noindent [Bi1] P. Biran, Symplectic packings in dimension 4,
Geom. and Funct. Anal. 7 (1997), no.3. 420-437.

\noindent [Bi2] P. Biran, Geometry of symplectic packing, PhD Thesis,
Tel-Aviv University (1997).

\noindent [Bi3] P. Biran, From symplectic packing to algebrain geometry
and back, to appear in the Proceedings of the 3'rd European Congress
of Mathematics.

\noindent [C] F. Catanese, Moduli space of surfaces
and real structures, math. AG/0103071.

%\noindent [D] S. Donaldson, Lefschetz fibrations in symplectic
%geometry, Doc. Math. J. DMV., ICM II (1998), 309-314.

\noindent [D] S. Donaldson, Symplectic submanifolds and almost-complex
geometry, J. Differential Geom. 44 (1996), no.4. 666-705.

\noindent [DS] S. Donaldson and I. Smith, in preparation.
%\item{} [Go] L. Gottsche, A conjectual generating function for
%numbers of curves on surfaces, preprint.

\noindent [FM1] R. Friedman and J. Morgan, On the diffeomorphism types
of certain algebraic surfaces,
J. Differential Geom. 27 (1988), no.3 371-398.

\noindent [FM2] R. Friedman and J. Morgan,
Algebraic surfaces and Seiberg-Witten invariants. J. Algebraic Geom. 6
(1997), no. 3, 445--479.

\noindent [FS1]  R. Fintushel and R. Stern,
Immersed spheres in 4-manifolds and the immersed Thom conjecture,
Turkish J. Math. 19 (1995), 145-157.

\noindent [FS2] R. Fintushel and R. Stern, Knots, links, and $4-$manifolds,
Invent. Math. 134 (1998) no. 2, 363-400.

\noindent [Ge] H. Geiges, Symplectic structures on $T^2-$bundles over $T^2$,
Duke Math. Journal, Vol. 67, no. 3 (1992), 539-555.

\noindent [Gr] M. Gromov, Pseudo-holomorphic curves in symplectic
manifolds, Invent. Math. 82 (1985), 307-347.

\noindent [IP] E-N. Ionel and T. Parker, The Gromov invariants of Ruan-Tian
and Taubes,  MRL. 4 (1997), 521-532.

\noindent [K] P. Kronheimer, Minimal genus in $S^1\times M^3$,
Invent. Math.
135 (1999) no.1 45-61.

\noindent [KK] V. Kharlamov and V. Kulikov, Real structures on RIGID surfaces, math.AG/0101098.

\noindent [La] F. Lalonde, Isotopy of symplectic balls,
Gromov's radius and the structure of ruled symplectic $4-$manifolds,
Math. Ann. 300 (1994), 273-296.

\noindent [Le] C. LeBrun, Diffeomorphisms, symplectic forms, and
Kodaira fibrations, preprint.

\noindent [L1] T. J. Li, Smoothly embedded spheres in symplectic
four manifolds,  Proc. Amer. Math. Soc. 127 (1999), no. 2, 609--613.

\noindent [L2] T. J. Li,
Symplectic Parshin-Arakelov inequality, Internat.
Math. Res. Notices 2000, no. 18, 941-954.

\noindent [Liu1] A. Liu,  Some new applications of
the general wall crossing formula, Math. Res. Letters 3 (1996), 569-585.

\noindent [Liu2] A. Liu, Family blow up formulas and nodal curves
in K\"ahler surfaces, preprint.

\noindent [LL1] T. J. Li and A. Liu, General wall crossing formula,
Math. Res. Letter 2 (1995), 797-810.

\noindent [LL2] T. J. Li and A. Liu, Symplectic
structures on ruled surfaces and a generalized adjunction inequality, Math.
Res. Letters 2 (1995), 453-471.

\noindent [LL3] T. J. Li and A. Liu, Family Seiberg-Witten invariants
and wall crossing formulas, to appear in Comm. in Analysis and Geometry.

\noindent [LL4] T. J. Li and A. Liu, On the equivalence between
SW and GT in the case $b^+=1$,
Internat.
Math. Res. Notices 1999, no. 7, 335--345.
%\item{} [LL5] T. J. Li and A. Liu, Counting curves on
%elliptic ruled surfaces, preprint.

\noindent [LiL1] B. H. Li and T. J. Li, Smooth minimal genera for
small negative classes in $CP^2\# n \overline{CP}^2$ when
$n\leq 9$, preprint.

\noindent [LiL2] B. H. Li and T. J. Li,
Symplectic genus, minimal genus and diffeomorphisms, preprint.

\noindent [LM] F. Lalonde and D. McDuff, The classification of ruled symplectic
$4-$manifolds, Math. Res. Lett. 3 (1996), no. 6, 769-778.

\noindent [Ma] M. Manetti, On the moduli space of diffeomorphic
algebraic surfaces, Inv. Math. (2001) to appear.

\noindent [M] J. Moser, On the volume elements on a manifold,
Trans. Amer. Math. Soc. 120 (1965), 286-294.

\noindent [Mc1] D. McDuff, The structure of rational and ruled
symplectic
$4$-manifolds. J. Amer. Math. Soc. 3 (1990), no. 3, 679--712.

\noindent [Mc2] D. McDuff, Lectures on
Gromov invariants for symplectic 4-manifolds,
With notes by Wladyslav Lorek. NATO Adv. Sci.
Inst. Ser. C Math. Phys. Sci., 488, Gauge theory and symplectic geometry
(Montreal, PQ, 1995), 175--210, Kluwer Acad. Publ., Dordrecht, 1997.

\noindent [Mc3] D. McDuff, From symplectic deformation to isotopy (Irvine, CA, 1996), 85-99, First Int. Press Lect Ser. I, Internat. Press, Cambridge, MA, 1998.

\noindent [Mc4] D. McDuff, Immersed spheres in symplectic $4-$manifolds,
Ann. Inst. Fourier (Grenoble) 42 (1992), no.1-2, 369-392.

\noindent [Mc5] D. McDuff, Blow ups and symplectic embeddings in dimension
4, Topology 30 (1991), 409-421.

\noindent [MP] D. McDuff and L. Polterovich, Symplectic packings and
algebraic geometry, Invent. Math. 115 (1994), no. 3. 405-434.

\noindent [MT] C. McMullen and  C. Taubes, $4-$manifolds with
inequivalent
symplectic forms and $3-$manifolds with inequivalent fibrations, preprint.

\noindent [R1] Y.B. Ruan, Symplectic topology on algebraic $3-$folds,
J. Differential Geom. 39 (1994), no. 1, 215-227.

\noindent [R2] Y.B. Ruan, Symplectic topology and extremal rays, Geom. Funct.
Anal. 3 (1993), no. 4, 395-430.

\noindent [R3] Y.B. Ruan, Symplectic topology and complex
surfaces. Geometry
and analysis on complex manifolds, 171--197, World Sci.
Publishing, River Edge, NJ, 1994.

\noindent [RT] Y.B. Ruan and G. Tian, A mathematical theory of
quantum
cohomology. J. Differential Geom. 42 (1995), no. 2, 259--367.

\noindent [S] D. Salamon, Spin geometry and Seiberg-Witten invariants,
book to appear.

\noindent [Sm] I. Smith, On moduli spaces of symplectic forms, preprint.

\noindent [T1] C.H. Taubes, {The Seiberg-Witten invariants and
   symplectic forms}, Math. Research Letters, 1 (1994) 809-822.

\noindent [T2] C.H. Taubes, {More constraints on symplectic
    manifolds from Seiberg--Witten invariants}, Math. Research
Letters, 2 (1995), 9-14.

\noindent [T3] C. Taubes, The Seiberg-Witten invariants and
the Gromov invariants, Math. Research Letters 2 (1995), 221-238.

\noindent [T4] C.H. Taubes, {Counting pseudo-holomorphic
submanifolds in dimension 4}, Journal of Diff. Geometry. 44 (1996) 818-893.

%\item{}[T2] C. H. Taubes, {\it SW$\Rightarrow$Gr: From Seiberg-Witten
%equations to pseudo-holomorphic curves}, JAMS, 9 (1996) 845-918.
%\item{}[T3] C. H. Taubes, {\it Gr$\Rightarrow$ SW: From pseudo-holomorphic
%curves to Seiberg-Witten solutions}, preprint.
%\item{}[T4] C. H. Taubes, {\it Gr=SW : Counting curves and connections}, prepr%int.

\noindent [Th] W. Thurston, Some simple examples of symplectic manifolds, Proc. Amer. Math. Soc. 55 (1976) 467-468.

\noindent [W] E. Witten, Monoples and Four-manifolds,
Math. Res. Letters 1
(1994) 769-796.

\noindent [Wi] P.M.H. Wilson, Symplectic deformations of Calabi-Yau
threefolds, J. Differential Geom. 45 (1997), no.3, 611-637.

\bigskip
\noindent Department of Mathematics, Princeton University, Princeton, NJ
08540

\noindent tli@math.princeton.edu

\medskip

\noindent Department of Mathematics, UC Berkeley, Berkeley, CA

\noindent akliu@math.berkeley.edu
\end